\documentclass[a4paper]{article}

\usepackage{a4wide}

\usepackage{amsmath,amssymb,latexsym,cite,amsthm,hyperref}

\usepackage{newtxtext}       %

\newcommand\bcmdtab{\noindent\bgroup\tabcolsep=0pt%
  \begin{tabular}{@{}p{10pc}@{}p{20pc}@{}}}
\newcommand\ecmdtab{\end{tabular}\egroup}

 \newtheorem{theorem}{Theorem}[section]
 \newtheorem{proposition}[theorem]{Proposition}
 \newtheorem{lemma}[theorem]{Lemma}
 \newtheorem{corollary}[theorem]{Corollary}
 \newtheorem{definition}[theorem]{Definition}

 \newtheorem{example}[theorem]{Example}

\usepackage{mathabx}
\usepackage[english]{babel}
\usepackage{tikz-cd}

\usepackage[shortcuts]{extdash}

\newcommand{\mat}[1]{{\ensuremath{\boldsymbol{\mathit{#1}}}}}
\newcommand{\alg}[1]{{\mathbf #1}}

\newcommand{\Mat}[1]{{\mathsf #1}}
\newcommand{\Fm}{{\ensuremath{\mathit{Fm}}}}

\newcommand{\dr}{\ensuremath{\textsc{dr}}}
\newcommand{\mdr}{\ensuremath{\textsc{mdr}}}
\newcommand{\acr}{\ensuremath{\textsc{acr}}}
\newcommand{\tcr}{\ensuremath{\textsc{tcr}}}
\newcommand{\gtcr}{\ensuremath{\textsc{gtcr}}}
\newcommand{\doo}{\ensuremath{\textsc{do}}}
\newcommand{\ds}{\ensuremath{\textsc{ds}}}

\newcommand\Mult{\text{\textit{\textbf{Mult}}}}
\newcommand\then{\DOTSB\;\Longrightarrow\;}
\newcommand\eqeq{\mathrel{\Dashv\mkern1.5mu\vDash}}
\newcommand\dedeq{\mathrel{\dashv\mkern1.5mu\vdash}}
\newcommand\FL{\mathcal{F\!L}}
\newcommand\Mod[1]{{#1\textnormal{-}\mathsf{Mod}}}
\newcommand*\cocolon{%
        \nobreak%
        \mskip6mu plus 1mu\mathpunct{}%
        \nonscript\mkern-\thinmuskip{:}\mskip2mu%
        \relax%
}
\let\models\vDash

\let\leq\leqslant
\let\geq\geqslant

\begin{document}

\begin{center}

{\Large An Abstract Approach to Consequence Relations}

\

\large

PETR CINTULA\\
            Institute of Computer Science of the Czech Academy of Sciences\\[1ex]
           JOS\'{E} GIL-F\'{E}REZ\\
           University of Bern\\[1ex]
            TOMMASO MORASCHINI\\
            Institute of Computer Science of the Czech Academy of Sciences\\[1ex]
             FRANCESCO PAOLI\\
            University of Cagliari

\end{center}

\noindent{\bf Abstract}\ {\small We generalise the Blok--J\'{o}nsson account of structural consequence relations, later developed by Galatos, Tsinakis and other authors, in such a way as to naturally accommodate multiset consequence. While Blok and J\'{o}nsson admit, in place of sheer formulas, a wider range of syntactic units to be manipulated in deductions (including sequents or equations), these objects are invariably \emph{aggregated} via set-theoretical union. Our approach is more general in that non-idempotent forms of premiss and conclusion aggregation, including multiset sum and fuzzy set union, are considered. In their abstract form, thus, \emph{deductive relations} are defined as additional compatible preorderings over certain partially ordered monoids. We investigate these relations using categorical methods and provide analogues of the main results obtained in the general theory of consequence relations. Then we focus on the driving example of multiset deductive relations, providing variations of the methods of matrix semantics and Hilbert systems in Abstract Algebraic Logic.}

\section{Introduction}

Most logicians agree that logic is about consequence, but different logicians have different preferences when it comes to specifying what this term means in their favourite propositional logics. Someone proceeds syntactically, taking consequence to be derivability in a given proof system, while someone else defines it semantically, using e.g.\ algebraic or Kripkean models. The most successful attempt to single out a neutral axiomatic framework embracing the common properties of all these different concepts is generally attributed to Tarski: in this acceptation, a consequence relation is a reflexive, monotonic and transitive relation between a set of formulas and a single formula of a given propositional language. Such Tarskian consequence relations (\tcr's: see Definition~\ref{def:TCR}) are the primary object of study of \emph{Abstract Algebraic Logic} (AAL: see e.g.~\cite{BP86,BP89,Cz01,Font,FoJa09,FJP03}), a discipline that aims at providing general tools for the investigation and comparison of the different brands of propositional logics on the market.

Among the reasons why AAL has established itself as a mainstream approach there is its success in effectively accommodating all the main extensions of, and alternatives to, classical propositional logic: intuitionistic logic, modal logics, relevant logics, quantum logics, and whatnot.

Although \tcr's are perfectly adequate for the needs of a wide spectrum of such abstract metalogical enquiries, it gradually emerged that they fail to capture a range of situations where we are still reasoning from given premisses to certain conclusions according to the same three principles of Reflexivity, Monotonicity, and Cut, yet we are not manipulating formulas of a given language, but perhaps \emph{sequents} (as in Gentzen calculi) or \emph{equations} (as in equational consequence relations associated with classes of algebras). Actually, it is not unusual for a logic to be given alternative presentations as a ``consequence relation'' of sorts over different sets of syntactic units. For example, classical logic can be presented not only as a (syntactically or semantically defined) \tcr, but also as the derivability relation of the classical sequent calculus, or as the equational consequence relation of the $2$-element Boolean algebra. In order to subsume these generalisations of the concept of propositional logic, core AAL was extended to $k$-dimensional systems (see~\cite{BP91}) and, subsequently, to Gentzen systems (see~\cite{Ra06,JRa13,ReV93,ReV95,Py99,T91}). The proliferation of these extensions of the classical AAL theory suggests that the idea of abstracting away from the specifics of these presentations to pinpoint what is essential to a logic is not without its allure.

\paragraph{Logical consequence: The Blok--J\'{o}nsson approach.}

In their groundbreaking paper~\cite{BJ} (see also the lecture notes~\cite{BlJo99} of the course ``Algebraic structures for logic'' that inspired the paper), Wim Blok and Bjarni J\'{o}nsson take their cue from such reflections and suggest to replace the set of formulas in the definition of \tcr\ by an \emph{arbitrary set} $A$, keeping however the same postulates of Reflexivity, Monotonicity, and Cut. These more general relations are called \emph{abstract consequence relations} (\acr's: see Definition~\ref{derkommissar}).

The first conceptual hurdle for this general approach is providing a suitable account of \emph{logical} consequence. It is generally agreed that logical consequence, as opposed to consequence in general, is a matter of logical form. In the standard AAL framework, this requirement is rendered precise by adding to the three Tarskian postulates the further condition of \emph{substitution\=/invariance}, also called sometimes \emph{structurality} (see Definition~\ref{def:TCR}). It is clear enough that replacing the algebra of formulas by the unstructured set $A$, which behaves as a sort of ``black box'' in so far as no relevant notion of endomorphism is applicable to it, calls for a completely different account of substitution\=/invariance.

Blok and J\'{o}nsson's response to this problem is insightful. They observe that the application of substitutions to propositional formulas (or, for that matter, to equations or sequents) is reminiscent of an operation of multiplication by a scalar. In fact, if $\varphi$ is a formula and $\sigma_{1},\sigma_{2}$ are substitutions, then $(\sigma_{1}\circ\sigma_{2})(\varphi) = \sigma_{1}(\sigma_{2}(\varphi))$ and, if $\iota$ is the identity substitution, $\iota(\varphi) = \varphi$. Generalising this example, we are led to an abstract counterpart of substitution invariance (see Definition~\ref{def:raclette}), where the set of substitutions is replaced by a monoid that acts on the set $A$.

One striking feature of Blok and J\'{o}nsson's suggestion is that it allows to reformulate, at a very general level, the classical notion of algebraisability in a way that provides a uniform perspective on the algebraisation of logics and Gentzen systems (see~\cite{BP89, ReV95}).

\paragraph{A categorical perspective.}

Nikolaos Galatos and Constantine Tsinakis (\cite{GT}) follow in Blok and J\'{o}nsson's footsteps and take their approach to the next level of generality. In particular, they aim at applying categorical and order-theoretic methods to the study of \acr's. One major hindrance to this accomplishment is the intrinsic \emph{asymmetry} of \acr's, whose relata are, respectively, a subset of the base set and an element thereof. With an eye to removing this potential source of technical problems, they introduce \emph{symmetric} versions of \acr's on a given set $A$. Thus, given a \acr\ $\vdash$ on a set $A$, they define its symmetric version $\vdash_s$ by $X\vdash_s Y$ iff $X\vdash y$ for all $y\in Y$. Although these symmetric \acr's, which coincide with those preorder relations on $\wp(A)$ containing the supersethood relation and such that $X \vdash \bigcup \{Z : X\vdash Z\}$, are shown to be in bijective correspondence with standard \acr's, their advantage is that their premiss\=/sets and conclusion\=/sets are points in a complete lattice of sets, $\wp(A)$. This circumstance suggests a natural generalisation to arbitrary complete lattices, here called \emph{Galatos--Tsinakis consequence relations} (\gtcr's: see Definition~\ref{def:gtcr}).

Can action-invariance be accommodated in this broader setting? To do so, we must first enrich our monoids of actions with additional structure, in such a way as to turn them into complete \emph{residuated lattices} (see~\cite{GJKO, MPT}). Monoidal actions are accordingly replaced by actions of such residuated lattices on complete lattices, giving rise to the notion of an $\mathbf M$-\emph{module} (see Definition~\ref{def:module}).

The algebraisability relation can be expressed in this framework in two different ways, respectively giving rise to two notions of equivalence between \gtcr's. One of the main results of the paper~\cite{GT} is the characterisation in purely categorical terms of the $\textbf M$\=/modules for which these two senses of equivalence coincide. These include the $\textbf M$\=/modules of formulas, of equations and of sequents. Thus, the result of Galatos and Tsinakis gives an elegant abstract characterisation of those $\textbf M$\=/modules that are just as ``well behaved'' as these standard examples.

\paragraph{The substructural challenge.}

As powerful and wide-ranging as this approach may be, a tricky challenge to its adequacy is posed by \emph{substructural logics} (see~\cite{Primer, GJKO, MPT}). Such logics, usually introduced by means of sequent calculi where some or all of the standard structural rules (Weakening, Contraction, Exchange, Cut) are restricted or even deleted, provide an interesting problem for the Blok--J\'{o}nsson approach to consequence, as further developed by Galatos and Tsinakis. By this we do not mean that they lie outside its scope --- on the contrary, they can be handled even in traditional AAL. There is, in fact, a canonical way to obtain a \tcr\ out of a given substructural sequent calculus. By way of example, consider the sequent calculus $\mathrm{FL}_{e}$ for \emph{full Lambek calculus with exchange} and the variety $\FL_{e}$ of pointed commutative residuated lattices. Define:
\begin{enumerate}
\item[$\bullet$] $X \vdash_{\mathrm{FL}_{e}} \varphi$ iff the sequent ${}\Rightarrow\varphi$ is provable in the calculus obtained by adding to $\mathrm{FL}_{e}$ as axioms the sequents in $\{{}\Rightarrow\psi:\psi\in X\}$ (this is sometimes called (see~\cite{Avronl}) the \emph{external consequence relation} of a sequent calculus).

\item[$\bullet$] $X\vdash_{\mathrm{FL}_{e}}' \varphi \quad\iff\quad \{\psi \wedge 1 \approx 1 : \psi\in X\} \models_{\FL_{e}}\varphi \wedge 1 \approx 1$.
\end{enumerate}
Here, $\models_\mathcal{K}$ refers to the \emph{equational consequence} of the class of algebras $\mathcal{K}$, namely, given a set of equations $\Pi$ and an equation $\varepsilon$, the expression $\Pi \models_\mathcal{K}\varepsilon$ means that every evaluation on an algebra of $\mathcal{K}$ satisfying all the equations of $\Pi$ also satisfies the equation $\varepsilon$.

It is well known \cite[Ch.~2]{GJKO} that both relations coincide, that is, ${\vdash_{\mathrm{FL}_{e}}} = {\vdash_{\mathrm{FL}_{e}}' }$. The relation defined by any of these two equivalent conditions is a \tcr\ to all intents and purposes.

This approach, however, seems to fly in the face of the motivation underlying substructural logics. \tcr's are relations between \emph{sets} of formulas and single formulas, whence they are insensitive to the number of occurrences a formula may have in some collection of premisses. In other words, they automatically validate the \emph{Contraction} and the \emph{Anticontraction} rules: if $X,\varphi,\varphi\vdash\psi$, then $X,\varphi\vdash\psi$, and conversely if $X,\varphi\vdash\psi$, then $X,\varphi,\varphi\vdash\psi$. Yet, some substructural logics (like \emph{linear logic},~\cite{Girard}) are commonly employed to formalise a ``resource-conscious'' notion of inference, according to which sentences are information tokens of a given type and for which the Contraction rule is utterly suspect. Other substructural logics in the relevant family (\cite{MMR, MMR2}) aim at capturing a concept of deduction according to which premisses in an argument should be actually \emph{used} to get the conclusion, something which seems to disqualify Anticontraction (and, even more, Weakening). In other words: AAL can certainly accommodate substructural logics into its framework (in the format of propositional logics or of Gentzen systems), and bestow on them the imposing bulk of general results it has to offer, but only at the cost of tweaking the substructural proof systems in such a way as to produce consequence relations that weaken and contract by \emph{fiat}. AAL, in sum, does not stay true to the spirit of substructural logics.

Now return, for a while, to the sequent calculus $\mathrm{FL}_{e}$. A more plausible candidate for a formalisation of substructural consequence is its so-called \emph{internal consequence relation} (see~\cite{Avronl}), namely, that relation that holds between a finite multiset of formulas $\Gamma$ and a formula $\varphi$ just in case $\Gamma\Rightarrow\varphi$ is a provable sequent of $\mathrm{FL}_{e}$. Investigating relations of this kind, however, means overstepping the Tarskian framework under at least two respects:

\begin{enumerate}
\item A consequence relation should be conceived of as a relation between a \emph{finite multiset} of formulas and a formula.

\item The Monotonicity postulate should be dropped and the Reflexivity postulate should be restricted.
\end{enumerate}

This approach, indeed, has been followed by Arnon Avron (see~\cite{Avron, Av92,Av94}) and, sporadically, by a few others (\cite{MMR, MMR2, Troelstra, Primer}), who laid down the fundamentals of a theory of multiset consequence. However, to help the theory to get started and make it easier to reconstruct some of the basic AAL theorems, it also seems wise to follow a middle-of-the-road perspective that shortens the gap with the Tarskian paradigm, adopting finite multisets as collections of premisses but leaving the Monotonicity and Reflexivity postulates untouched. The resulting relations have a built-in Weakening condition, although they do not necessarily contract. This policy, as a matter of fact, faces an insurmountable problem. David Ripley (see~\cite{Ripley}) has shown that it is not possible to obtain a bijective correspondence between these ``naive'' multiset consequence relations and closure operators on finite multisets of formulas --- any such relation that arises from a closure operator has to obey Contraction.

In the paper~\cite{CP}, the authors developed a strategy to avoid this problem. They adopted a multiple-conclusion format, studying relations $\vdash$ between finite multisets of formulas of a given propositional language, here called \emph{multiset deductive relations} (\mdr's: see Definition~\ref{d:MDR}). As in the case of \gtcr's, the interpretation of the right-hand side is essentially conjunctive (cf.\ similar approach in non-monotonic case in~\cite{Av92}). In other words, $\Gamma \vdash \Delta$ can be read as: using at most once all the formula occurrences in $\Gamma$, we can derive all the formula occurrences in $\Delta$. The paper contains arguments in favour of this preference over a disjunctive reading. They also suitably modified the notions of a closure operator and a closure system on finite multisets of formulas so as to recover the traditional lattice isomorphism results that characterise the standard set\=/theoretical framework. These correspondences were laid down as the embryo of a theory that aimed at eventually obtaining appropriate analogues of the main results available in AAL for Tarskian consequence.

It is worth asking whether something like these multiset consequence relations, which have no direct counterpart in AAL, can still ensconce in some more flexible apparatus based on (not necessarily complete) lattices. After all, finite multisets over a set still form a lattice under the operations
\[
 (\mathfrak{X}\vee\mathfrak{Y})(a) = \sup\{\mathfrak{X}(a),\mathfrak{Y}(a)\}
 \quad\text{and}\quad
 (\mathfrak{X}\wedge\mathfrak{Y})(a) = \inf\{\mathfrak{X}(a),\mathfrak{Y}(a)\}
\]
(see below), out of which we could define an $\mathbf{M}$\=/module of sorts that would land us in known territory. However, it is not hard to see that the crucial operation on multisets is not any of these, but rather \emph{multiset sum}:
\[
 (\mathfrak{X}\uplus\mathfrak{Y})(a) = \mathfrak{X}(a)+\mathfrak{Y}(a).
\]
It is via multiset sum that we aggregate multisets of premisses in substructural logics and formulate sequent rules in substructural sequent calculi. Being a non-idempotent operation, though, it is scarcely pliant to the methods reviewed so far. The \gtcr's on complete lattices, therefore, need to be replaced by appropriate relations on \emph{dually integral partially ordered Abelian monoids}, the prime motivating example being the pomonoid
\[
 \langle \Fm_{\mathcal{L}}^{\flat}, \leqslant, \uplus,\emptyset \rangle,
\]
where $\Fm_{\mathcal{L}}^{\flat}$ is the set of finite multisets of formulas of a propositional language $\mathcal{L}$, $\uplus$ is multiset sum, $\emptyset$ is the empty multiset and $\mathcal{\leqslant}$ is the sub\=/multisethood relation (all these notions will be rehearsed in Subsection~\ref{carmela}).

The study of such \emph{deductive relations} (\dr's: see Definition~\ref{Def:DeductiveRelation}) is the main topic of the present paper. Before outlining its contents in Subsection~\ref{rosalia}, however, it is expedient to go through a number of preliminaries needed to make the paper reasonably self-contained.

\section{Preliminaries}

We start this section by reviewing some basic notions about multisets only to such an extent as it is needed for the purposes of the present paper. For a more comprehensive account, the reader can consult e.g.~\cite{Blizard, SYS08}.

\subsection{Multisets}\label{carmela}

By a \emph{multiset} over a set $A$ we mean a function $\mathfrak{X}$ from $A$ to the set $\mathbb{N}$ of natural numbers.\footnote{We use letters $\Gamma,\Delta,\Pi,\dots$ for multisets of \emph{formulas}, while $\mathfrak{X}$, $\mathfrak{Y}$, $\mathfrak{Z},\dots$ are used for general multisets.} By $\wp^{M}(A)$ we denote the set of all multisets over $A$. The \emph{root set} of a multiset $\mathfrak{X}$ is the set
\[
 \vert \mathfrak{X} \vert = \{a\in A : \mathfrak{X}(a) > 0\}.
\]
If $a\in \vert \mathfrak{X} \vert$, we say that $a$ is an \emph{element} of $\mathfrak{X}$ of \emph{multiplicity} $\mathfrak{X}(a)$. A multiset $\mathfrak{X}$ is \emph{finite} if $\vert\mathfrak{X}\vert$ is finite. By $A^{\flat}$ we denote the set of all finite multisets over $A$. The \emph{empty multiset}, i.e.\ the constant function $0$, will be denoted by the same symbol $\emptyset$ used for the empty set --- the context will always be sufficient to resolve ambiguities.

The set $\wp^{M}(A)$ inherits the ordering of $\mathbb{N}$ in the following way:
\[
 \mathfrak{Y} \leqslant \mathfrak{X} \quad\iff\quad
 \mathfrak{Y}(a) \leq \mathfrak{X}(a), \text{ for all } a\in A.
\]
With respect to this ordering, it forms a lattice with joins and meets defined as
\[
 (\mathfrak{X}\vee\mathfrak{Y})(a) = \sup\{\mathfrak{X}(a),\mathfrak{Y}(a)\}
 \quad\text{and}\quad
 (\mathfrak{X}\wedge\mathfrak{Y})(a) = \inf\{\mathfrak{X}(a),\mathfrak{Y}(a)\},
\]
for all $a\in A$. The operation $\vee$ is a kind of ``union'', and, true to form, if we consider the subsets of $A$ as multisets whose elements have multiplicity $1$ and $\mathfrak{X}$ and $\mathfrak{Y}$ are subsets of $A$, then $\mathfrak{X}\vee\mathfrak{Y} = \mathfrak{X}\cup\mathfrak{Y}$. There is another ``union\=/like'' operation of \emph{sum} between multisets, defined as follows:
\[
 (\mathfrak{X}\uplus\mathfrak{Y})(a)=\mathfrak{X}(a)+\mathfrak{Y}(a),\quad \text{for all } a\in A.
\]

The next proposition shows that the set of finite multisets over a set $A$ can be seen as the universe of a dually integral Abelian pomonoid which we will denote as $\mathbf{A}^{\flat}$.

\begin{proposition}\label{p:MutliPomon}
For any set $A$ the structure
$
 \mathbf{A}^{\flat} = \langle A^{\flat},\leqslant,\uplus,\emptyset \rangle
$
is a dually integral Abelian pomonoid, i.e., a structure where:
\begin{enumerate}
\item $\langle A^{\flat},\uplus,\emptyset\rangle$ is an Abelian monoid.
\item $\leqslant$ is a partial order compatible with $\uplus$, i.e.,
\[
 \mathfrak{X}\leqslant\mathfrak{Y} \quad\then\quad
 \mathfrak{X}\uplus\mathfrak{Z}\leqslant \mathfrak{Y}\uplus\mathfrak{Z}.
\]
\item $\emptyset$ is the bottom element of\/ $\leqslant$.
\end{enumerate}
\end{proposition}

We will also have occasion to use the operation $\mathfrak{X}\setminus\mathfrak{Y}$ defined by
\[
 (\mathfrak{X}\setminus\mathfrak{Y})(a) =
 \max\bigl\{\mathfrak{X}(a)-\mathfrak{Y}(a), 0 \bigr\},\quad \text{for all } a\in A,
\]
relying on the context to disambiguate between this operation and standard set-theoretic subtraction. As it is customary to do, we use square brackets for multiset abstraction; so, for example, $[a,a,b,c]$ will denote the multiset $\mathfrak{X}$ s.t.\ $\mathfrak{X}(a)=2$, $\mathfrak{X}(b)=\mathfrak{X}(c)=1$, and $\mathfrak{X}(d)=0$, for any $d\notin\{a,b,c\}$.

Finally, every map $f\colon A\to B$ can be extended to a morphism from $\mathbf{A}^\flat$ to $\mathbf{B}^\flat$ (for which we retain the same symbol) via
\[
  f(\mathfrak{X}) = [f(a_{1}), \dots, f(a_{n})]
\]
for every $\mathfrak{X} = [a_{1}, \dots, a_{n}]$. We will use this notation without special mention.

\subsection{Consequence relations}

Let $\mathcal{L}$ be a propositional language (or, which is the same, an algebraic language). By $\mathbf{Fm}_{\mathcal L}$ we denote the algebra of formulas with countably many variables on the language $\mathcal L$ with universe $Fm_{\mathcal L}$. An $\mathcal{L}$-\emph{substitution} is an endomorphism $\sigma$ of $\mathbf{Fm}_{\mathcal{L}}$. The notation $\sigma(X)$ refers to the set $\{\sigma(\varphi) : \varphi\in X\}$, for every $X\subseteq Fm_{\mathcal{L}}$.

\begin{definition}\label{def:TCR}
A \emph{Tarskian consequence relation} (\tcr) on $\mathcal{L}$ is a binary relation ${\vdash} \subseteq \wp(\Fm_{\mathcal{L}})\times \Fm_{\mathcal{L}}$ obeying the following conditions for all $X,Y\subseteq \Fm_{\mathcal{L}}$ and $\varphi\in\Fm_{\mathcal{L}}$:
\begin{enumerate}
\item[$\bullet$] $X\vdash \varphi$ whenever $\varphi\in X$. \hfill (\emph{Reflexivity})

\item[$\bullet$] If\/ $X\vdash \varphi$ and $X\subseteq Y$, then $Y\vdash \varphi$. \hfill (\emph{Monotonicity})

\item[$\bullet$] If\/ $Y\vdash \varphi$ and $X\vdash \psi$ for every $\psi\in Y$, then $X\vdash \varphi$. \hfill (\emph{Cut})
\end{enumerate}
A \tcr\ $\vdash$ on $\mathcal{L}$ is said to be \emph{substitution\=/invariant}, or also a \emph{logic}, if for all $X\cup\{\varphi\} \subseteq \Fm_{\mathcal{L}}$, whenever $X\vdash\varphi$ we also have that $\sigma(X) \vdash \sigma(\varphi)$, where $\sigma$ is an arbitrary $\mathcal{L}$-substitution.
\end{definition}

As we mentioned in the introduction, \tcr's are the main object of study of AAL. The next definition, due to Blok and J\'{o}nsson (\cite{BJ}) generalises \tcr's so as to also encompass cases in which the syntactic units are not formulas, e.g.\ equational consequence relations and consequence relations on sequents.

\begin{definition}\label{derkommissar}
An \emph{abstract consequence relation} (\acr) on a set $A$ is a relation ${\vdash}\subseteq{\wp(A)\times A}$ obeying the following conditions for all $X,Y\subseteq A$ and for all $a\in A$:
\begin{enumerate}
\item[$\bullet$] $X\vdash a$ whenever $a\in X$. \hfill (\emph{Reflexivity})

\item[$\bullet$] If\/ $X\vdash a$ and $X\subseteq Y$, then $Y\vdash a$. \hfill (\emph{Monotonicity})

\item[$\bullet$] If\/ $Y\vdash a$ and $X\vdash b$ for every $b\in Y$, then $X\vdash a$. \hfill (\emph{Cut})
\end{enumerate}
\end{definition}

In this framework, generalising substitution invariance requires a little more work. The set of $\mathcal L$-substitutions, which forms a monoid with composition, is replaced by an arbitrary monoidal action on the base set $A$.

\begin{definition}\label{def:raclette}
A monoid $\mathbf{M}=\langle M,\cdot,1\rangle$ \emph{acts} on a set $A$ if there is a map $\star\colon M\times A \to A$ s.t.\ for all $m_{1},m_{2}\in M$ and all $a\in A$,
\[
 (m_{1}\cdot m_{2})  \star a=m_{1}\star(m_{2}\star a)  \quad\text{and}\quad 1\star a=a.
\]
Thus, an \acr\ $\vdash$ is said to be \emph{action-invariant} if for all $X\cup\{a\}  \subseteq A$ and $m\in M$, whenever $X\vdash a$ we also have that $\{m\star x : x\in X\}  \vdash m\star a$.
\end{definition}

Observe that if $\mathbf{M}$ acts on $A$, then $\wp(\mathbf{M}) = \langle \wp(M),\cdot' ,\{1\}\rangle$ (where $\cdot'$ is complex product) acts on $\wp(A)$ via the induced map
\[
  N\star' X = \{m\star x : m\in N,\ x\in X\}.
\]
The sets $\wp(A)$ and $\wp(\mathbf M)$ have the structures of a complete lattice and a complete residuated lattice under set\=/inclusion, respectively. Moreover, the map $\star' \colon \wp(M)\times \wp(A)\to \wp(A)$ is \emph{biresiduated}, i.e., it is residuated in each coordinate. This motivates the following definitions.

\begin{definition}\label{def:module}
Let $\mathbf{M}=\langle M,\wedge^{\mathbf{M}},\vee^{\mathbf{M}},\cdot,\backslash,/,1\rangle $ be a complete residuated
lattice, $\mathbf{L}=\langle L,\wedge^{\mathbf{L}},\vee^{\mathbf{L}}\rangle $ be a complete lattice and $\star\colon M\times L\to L$ be a map. We say that $\mathbf{M}$ \emph{acts on} $\mathbf{L}$, or also that $\boldsymbol{L} = \langle L,\wedge^{\mathbf{L}},\vee^{\mathbf{L}},\star\rangle$ is an $\mathbf{M}$-\emph{module},\footnote{Note that for modules we use boldface italic font, whereas for algebras a simple boldface. Also, if $\boldsymbol{K}$ is module then $\mathbf{K}$ is its lattice reduct.} if the monoid reduct of $\mathbf{M}$ acts on $L$ and, moreover, there are maps $/_{\star}\colon L\times L\to M$ and $\backslash_{\star}\colon M\times L\to L$ such that, for all $m\in M$ and $x,y\in L$,
\[
 m\star x\leq^{\boldsymbol{L}}y \quad\iff\quad
 m \leq^{\mathbf{M}}y/_{\star}x \quad\iff\quad
 x\leq^{\boldsymbol{L}}m\backslash_{\star}y.
\]
\end{definition}

Every \acr\ $\vdash$ on a set $A$ can be lifted to a binary relation $\vdash_s$ on $\wp(A)$ as follows:
\[
  X\vdash_s Y \quad\iff\quad X\vdash y, \text{for every } y\in Y.
\]
Observe that $\vdash_s$ satisfies these three properties:
\begin{enumerate}
\item[$\bullet$] $\vdash_s$ is a preorder on $\wp(A)$,
\item[$\bullet$] if $X\subseteq Y$ then $Y\vdash_s X$,
\item[$\bullet$] $X\vdash\bigcup\{Z : X\vdash Z\}$.
\end{enumerate}
Also, any binary relation on $\wp(A)$ satisfying these three properties is of the form $\vdash_s$ for some \acr\ $\vdash$ on $A$, and this correspondence is bijective. If moreover $\vdash$ is an action-invariant \acr, then $\vdash_s$ satisfies that
\[
 X\vdash_s Y \quad\then\quad N\star' X \vdash_s N\star' Y.
\]
It is natural enough to propose an abstract version of these relations, as done in~\cite{GT}.

\begin{definition}\label{def:gtcr}
A \emph{Galatos--Tsinakis consequence relation} (\gtcr) on a complete lattice $\mathbf{L} = \langle L,\wedge,\vee\rangle $, with induced order $\leq$, is a preorder $\vdash$ of $L$ that contains $\geq$ and is such that for all $x\in L$, $x\vdash\bigvee\{y\in L : x\vdash y\}$. A \gtcr\ $\vdash$ on an $\mathbf{M}$\=/module $\boldsymbol{L}$ is said to be \emph{action-invariant}  if, for all $x,y\in L$ and $m\in M$, $x\vdash y$ implies $m\star x\vdash m\star y$.
\end{definition}

In traditional AAL, the fundamental concept of deductive equivalence is the algebraisability relation (\cite{BP89,Font}) between logics and equational consequence relations. We recall its definition hereafter.

\begin{definition}\label{def:fondue}
A logic $\vdash$ on a language $\mathcal L$ is \emph{algebraisable} if there exist a generalised quasi\=/variety $\mathcal{K}$ and two maps
\[
 \tau\colon\wp(\Fm_{\mathcal{L}}) \longleftrightarrow \wp(Eq_{\mathcal{L}})\cocolon\rho
\]
which commute with unions and substitutions, such that
\[
 X\vdash\varphi \iff \tau(X)\models_{\mathcal{K}}\tau(\varphi)
 \quad\text{and}\quad
x\thickapprox y \eqeq_{\mathcal K} \tau\rho(x\thickapprox y)
\]
for every $X\cup\{\varphi\}\subseteq \Fm_{\mathcal{L}}$.
\end{definition}

 In the displays above, $Eq_{\mathcal{L}}$ is the set of $\mathcal{L}$\=/equations (formally cast as ordered pairs of $\mathcal{L}$\=/formulas), and~$\models_{\mathcal{K}}$ the equational consequence relation of $\mathcal{K}$. When the above conditions hold, $\mathcal{K}$ is unique and is said to be the \emph{equivalent algebraic semantics} of $\vdash$. Given a logic $\vdash$ and a generalised quasi-variety $\mathcal{K}$, we denote by $C_{\vdash}$ and $C_{\mathcal{K}}$ the closure operators associated to the relations $\vdash$ and $\models_{\mathcal{K}}$ respectively. It is well known (see e.g.~\cite[p.~149]{Font}) that the notion of algebraisability is captured by the existence of a particular isomorphism. Below, given an \acr\ $\vdash$ on a set $A$, we denote by $\mathcal{T}h(\vdash)$ the complete lattice of \emph{theories} of $\vdash$, i.e., such sets $T\subseteq A$ such that $a\in T$ whenever $T\vdash a$.

\begin{theorem}[Syntactic Isomorphism Theorem]\label{Thm:SyntacticIsoThm}
Let\/ $\vdash$ be a logic and\/ $\mathcal{K}$ a generalised quasi\=/variety. Then\/ $\vdash$ is algebraisable with equivalent algebraic semantics\/ $\mathcal{K}$ if and only if there is a lattice isomorphism\/ $\Phi\colon \mathcal Th(\vdash)\to \mathcal Th(\models_{\mathcal{K}})$ such that\/ $\Phi\circ C_{\vdash}\circ\sigma=C_{\mathcal{K}}\circ\sigma\circ\Phi$ for every substitution $\sigma$.
\end{theorem}

Blok and J\'{o}nsson introduced the following natural generalisation of the concept of algebraisability \cite[Definition~4.5]{BlJo99}, based on the criterion provided by Theorem~\ref{Thm:SyntacticIsoThm}.

\begin{definition}
Let $\mathbf{M}$ be a monoid acting on two sets $A_{1}$ and $A_{2}$, respectively through the actions $\star_{1}$ and $\star_{2}$. Moreover, let $\vdash_{1}$ and $\vdash_{2}$ be two action invariant \acr's, respectively on $A_{1}$ and $A_{2}$. The \acr's $\vdash_{1}$ and $\vdash_{2}$ are said to be \emph{equivalent} if the following two lattices with unary operators are isomorphic
\[
 \langle \mathcal Th(\vdash_1), C_{\vdash_1} m  : m \in M\rangle \cong
 \langle \mathcal Th(\vdash_2), C_{\vdash_2} m  : m \in M\rangle.
\]
\end{definition}

Observe that, in the light of Theorem~\ref{Thm:SyntacticIsoThm}, a logic $\vdash$ is algebraisable with equivalent algebraic semantics $\mathcal{K}$ if and only if the substitution\=/invariant \acr's $\vdash$ and $\models_{\mathcal{K}}$ are equivalent according to this definition (where $\mathbf{M} = \mathbf{End}(\mathbf{Fm}_{\mathcal{L}}) $ is the monoid of substitutions, acting in the natural way on formulas and equations).

Since the maps $\tau\colon\wp(\Fm_{\mathcal{L}})\longleftrightarrow\wp(Eq_{\mathcal{L}})\cocolon\rho$ in Definition~\ref{def:fondue} commute with arbitrary unions, they are residuated maps on the corresponding complete lattices. Moreover, they commute with substitutions. This is the reason why \cite{GT} consider the category $\Mod{\mathbf M}$ whose objects are $\mathbf{M}$\=/modules and whose arrows are residuated maps $\tau\colon\boldsymbol{L}\to \boldsymbol{L}' $ (called \emph{translators}) such that, if $\mathbf{M}$ acts on $\mathbf{L}$ and $\mathbf{L}' $ via $\star_{1}$ and $\star_{2}$ respectively, we have that for all $x\in L$ and $m\in M$, $\tau(m\star_{1} x) = m\star_{2}\tau(x)$. A noteworthy feature of $\Mod{\mathbf{M}}$ is the fact that \gtcr's on its objects can be viewed as \emph{bona fide} objects in the same category.

\begin{theorem}\mbox{}
\begin{enumerate}
\item The \gtcr's on an\/ $\mathbf{M}$\=/module $\boldsymbol{L}$ correspond bijectively to closure operators on its lattice reduct\/ $\mathbf{L}$ (namely, enlarging, order\=/preserving, and idempotent unary operation\/ $\gamma$ on\/ $\mathbf{L}$) via the maps\/ $\vdash_{()}$ and\/ $\gamma_{()}$ defined by
\[
 x\vdash_{\gamma}y\text{ iff } y\leq^{\mathbf{L}}\gamma(x)  \quad\text{and}\quad
 \gamma_{\vdash}(x)  =\bigvee\{y\in L : x\vdash y\}.%
\]

\item If $\boldsymbol{L}$ is an\/ $\mathbf{M}$\=/module (via $\star$) and\/ $\vdash$ is a \gtcr\ on its lattice reduct\/ $\mathbf{L}$, then the lattice\/ $\mathbf{L}_{\gamma_{\vdash}}$ of\/ $\gamma_{\vdash}$-closed elements of\/ $\mathbf{L}$ is the lattice reduct of an\/ $\mathbf{M}$\=/module $\boldsymbol{L}_{\gamma_{\vdash}}$ via the map\/ $\star_{\gamma_{\vdash}}\colon M\times L_{\gamma_{\vdash}}\to L_{\gamma_{\vdash}}$ defined by
\[
 m\star_{\gamma_{\vdash}}x = \gamma_{\vdash}(m\star x).
\]

Moreover, $\gamma_{\vdash}$ is a morphism in\/ $\Mod{\mathbf{M}}$ from $\boldsymbol{L}$ onto $\boldsymbol{L}_{\gamma_{\vdash}}$.
\end{enumerate}
\end{theorem}

Observe that if $\vdash$ is a \gtcr\ on $\boldsymbol{L}$, then $\mathbf{L}_{\gamma_{\vdash}}$ is nothing but the lattice of $\vdash$\=/\emph{theories} --- viz., of all $t\in L$ such that $t\vdash x$ implies $x\leq t$.

It turns out that algebraisability can be generalised to the setting of modules over residuated lattices, again thanks to the criterion provided by Theorem~\ref{Thm:SyntacticIsoThm}.

\begin{definition}
Two \gtcr's $\vdash_{1}$ and $\vdash_{2}$, respectively on the $\mathbf{M}$\=/modules $\boldsymbol{L}$ and $\boldsymbol{L}' $, are said to be \emph{equivalent} if there is a module isomorphism $f\colon \boldsymbol{L}_{\gamma_{\vdash_{1}}}\to\boldsymbol{L}'_{\gamma_{\vdash_{2}}}$. An isomorphism $f\colon\boldsymbol{L}_{\gamma_{\vdash_{1}}}\to\boldsymbol{L}'_{\gamma_{\vdash_{2}}}$ is \emph{induced} by the translators $\tau\colon\boldsymbol{L}\to\boldsymbol{L}' $ and $\rho\colon\boldsymbol{L}' \to \boldsymbol{L}$ if $f\gamma_{\vdash_{1}}=\gamma_{\vdash_{2}}\tau$ and $f^{-1}\gamma_{\vdash_{2}}=\gamma_{\vdash_{1}}\rho$.
\end{definition}

Given an isomorphism $f\colon\boldsymbol{L}_{\gamma_{\vdash_{1}}}\to\boldsymbol{L}'_{\gamma_{\vdash_{2}}}$ induced by translators $\tau$ and $\rho$, the classical definition of algebraisability can be restored, in the sense that for every $x,y\in L$ and $z\in L' $ we have that
\[
 x\vdash_{1}y \iff \tau(x) \vdash_{2} \tau(y)
 \quad\text{and}\quad
 z \dedeq_{2}\tau(\rho(z)).
\]

In this parlance, the Syntactic Isomorphism Theorem states that every equivalence between a logic $\vdash$ and the equational consequence $\models_{\mathcal{K}}$ relative to a generalised quasi-variety $\mathcal{K}$ is induced by a pair of translators. It is natural to ask whether this is true for arbitrary equivalences between \gtcr's on $\mathbf{M}$\=/modules. Unfortunately, it turns out that this is false in general, as shown in~\cite{GilFe11}. Nevertheless, Galatos and Tsinakis find sufficient and necessary conditions for it to be the case. Recall, that an object $\boldsymbol{R}$ in a category $\mathsf{C}$, whose arrows are set-theoretic functions, is \emph{onto\=/projective} if for any $\mathsf{C}$-morphisms $f\colon\boldsymbol{S} \to\boldsymbol{T}$ and $g\colon\boldsymbol{R}\to \boldsymbol{T}$ between objects in $\mathsf{C}$ with $f$ onto, there is a $\mathsf{C}$-morphism $h\colon\boldsymbol{R}\to\boldsymbol{S}$ such that $f\circ h=g$. Here is the main result in the paper by Galatos and Tsinakis:

\begin{theorem}\label{calabro}
An $\mathbf{M}$\=/module $\boldsymbol{L}$ is onto\=/projective in $\Mod{\mathbf{M}}$ if and only if for any $\mathbf{M}$\=/module $\boldsymbol{L}' $ and \gtcr's $\vdash_{1}$ and $\vdash_{2}$, respectively on $\boldsymbol{L}$ and $\boldsymbol{L}' $, every residuated order embedding $f \colon \boldsymbol{L}_{\gamma_{\vdash_{1}}} \to \boldsymbol{L}'_{\gamma_{\vdash_{2}}}$ is such that $f\gamma_{\vdash_{1}} = \gamma_{\vdash_{2}}\tau$ for some translator $\tau\colon \boldsymbol{L}\to \boldsymbol{L}'$. In particular, if both $\boldsymbol{L}$ and $\boldsymbol{L}'$ are onto\=/projective, then every equivalence between $\vdash_{1}$ and $\vdash_{2}$ is induced by translators.
\end{theorem}

\noindent Generalisations of the above criterion have been obtained in~\cite{GaGF1x-a} (see also~\cite{CRu13}) and~\cite{Mor16}, while the structure of onto\=/projective objects in the Blok--J\'{o}nsson framework was described in~\cite{FM13}. Note that in~\cite{GaGF1x-a} the authors have also proved that all epimorphisms in $\Mod{\mathbf{M}}$ are onto, so the specification ``onto'' in the previous theorem is redundant.

Given a propositional language $\mathcal{L}$, the lattice of sets $\wp(\mathbf{Fm}_{\mathcal{L}})$ corresponding to the plain old module of $\mathcal{L}$\=/formulas is indeed onto\=/projective in $\Mod{\wp(\mathbf{End}(\mathbf{Fm}_{\mathcal{L}}))}$, and so are the modules of $\mathcal{L}$\=/equations and $\mathcal{L}$\=/sequents. Thus, Theorem~\ref{calabro} identifies modules that are just as ``well behaved'' as these standard examples in terms of admitting a general version of the Syntactic Isomorphism Theorem.

The following concept, already discussed in the introduction, is the key motivating example for the further generalisations that will be at the centre of the next sections.

\begin{definition}\label{d:MDR}
A \emph{multiset deductive relation} (\mdr) on a propositional language $\mathcal{L}$ is a relation $\vdash$ on $\Fm_{\mathcal{L}}^{\flat}$ such that for each $\Gamma,\Delta,\Pi\in \Fm_{\mathcal{L}}^{\flat}$:
\begin{enumerate}
\item[$\bullet$] $\Gamma\uplus\Delta\vdash \Gamma$. \hfill (\emph{Reflexivity})

\item[$\bullet$] If\/ $\Gamma\vdash \Delta$, then $\Gamma\uplus\Pi\vdash \Delta\uplus\Pi$. \hfill (\emph{Compatibility})

\item[$\bullet$] If\/ $\Gamma\vdash \Delta$ and $\Delta\vdash \Pi$, then $\Gamma\vdash \Pi$. \hfill (\emph{Transitivity})
\end{enumerate}
\end{definition}

\subsection{Overview of the results}\label{rosalia}

This paper is structured as follows. In Section~\ref{dudu}, we introduce the concept of a \emph{deductive relation} (\dr) on a dually integral Abelian pomonoid. This is a modification of the notion of \gtcr\ on a complete lattice, so as to encompass \mdr's and other examples that are not directly covered by Galatos and Tsinakis' theory. Since, as we have seen, \gtcr's on a lattice are in bijective correspondence with closure operators on the same lattice, it is to be expected --- if we are on the right track --- that an analogous result holds with respect to \dr's and some sort of ``operational companions'' of such. The fact that aggregation of premisses and conclusions is abstractly represented by a monoidal operation which, unlike set union, is not necessarily idempotent, implies that it won't do to define these operators in the standard way. Thus, \emph{deductive operators} (\doo's) on a dually integral Abelian pomonoid $\mathbf{R} = \langle R,\leq,+,0\rangle$ are introduced as certain maps $\delta\colon R\to \wp(R)  $. Similarly, we propose a notion of a \emph{deductive system} (\ds) that appropriately generalises closure systems associated with closure operators. The main result of the section is:

\begin{theorem}[see Theorem~\ref{Thm:Trinity}]
If\/ $\mathbf{R}$ is a dually integral Abelian pomonoid, then the posets\/ $\langle \mathrm{\operatorname{Rel}}(\mathbf{R}),\subseteq \rangle$, $\langle \mathrm{Op}(\mathbf{R}),\preccurlyeq\rangle $ and\/ $\langle \mathrm{Sys}(\mathbf{R}),\supseteq\rangle$ are isomorphic.
\end{theorem}

In Section~\ref{azione}, the problem of action-invariance is under scrutiny. We define a category $\Mod{\mathbf{A}}$ of modules over dually integral po\=/semirings (called $\mathbf{A}$-\emph{modules}) that is closely related to the category of Galatos and Tsinakis' $\mathbf{M}$\=/modules and includes as new examples the modules $\Mult_{\mathcal{L}}$ whose underlying pomonoids have the form $\langle \Fm_{\mathcal{L}}^{\flat},\mathcal{\leqslant},\uplus,\emptyset\rangle $ for some language $\mathcal{L}$, and whose scalars are finite multisets of $\mathcal{L}$-substitutions. In this wider framework, we obtain analogues of the main theorems proved by Galatos and Tsinakis. Here is an example:

\begin{theorem}[see Theorem~\ref{nacchera}]
An $\mathbf{A}$\=/module $\boldsymbol{R}$ is onto\=/projective in  $\Mod{\mathbf{A}}$ iff, for any other $\mathbf{A}$\=/module $\mathbf{S}$ and action-invariant \doo's $\delta$ and $\gamma$ on $\mathbf{R}$ and $\mathbf{S}$ respectively, every injective and order-reflecting morphism $\Phi\colon\mathbf{R}_{\delta} \to \mathbf{S}_{\gamma}$ is induced by some morphism.
\end{theorem}

In Theorem~B, $\mathbf{R}_{\delta},\mathbf{S}_{\gamma}$ are $\mathbf{A}$\=/modules of sets whose universes are, respectively, the union of all $\delta$-images (resp., $\gamma$-images) of elements of $R$ (resp.,~$S$). We also show that our motivating multiset-theoretical example is just as well behaved as the standard examples in Galatos and Tsinakis' theory:

\begin{theorem}[see Theorem~\ref{blascoferrer}]
Modules arising from multisets are onto\=/projective in the appropriate categories.
\end{theorem}

In Section~\ref{energade}, we zoom in on \mdr's. We introduce two types of matrix semantics for such relations. If $\vdash$ is an \mdr\ on a language $\mathcal{L}$, an $\mathcal{L}$-\emph{hypermatrix} is a pair $\langle \mathbf{A},F\rangle $, where $\mathbf{A}$ is an $\mathcal{L}$-algebra and $F$ is a $\leqslant$-downward closed set of finite multisets over $A$. We show that the $\mathcal{L} $-hypermatrix models of an \mdr\ $\vdash$ correspond bijectively to the models of a Gentzen relation uniquely associated to $\vdash$. This opens the way for importing into our theory all sorts of tools and results from the abstract theory of Gentzen systems (see~\cite{Ra06, ReV93, ReV95}), including a workable definition of Leibniz congruence of an $\mathcal{L}$-hypermatrix and a completeness theorem for any substitution\=/invariant \mdr:

\begin{theorem}[see Theorem~\ref{t:CompletenessReduced}]
Every substitution\=/invariant \mdr\ is complete with respect to the class of its reduced $\mathcal{L}$-hypermatrix models.
\end{theorem}

An alternative type of matrix semantics for \mdr's, on the other hand, involves certain structures made up by an algebra $\mathbf{A}$, a dually integral Abelian pomonoid $\langle D,\leq,+,0\rangle $,\textbf{ }and a pomonoid homomorphism from $\langle A^{\flat},\mathcal{\leqslant},\uplus,\emptyset\rangle $ to $\langle D,\leq,+,0\rangle $. These structures subsume ordinary logical matrices, which arise when $\mathbf{D}$ is the $2$-element join semilattice. We clarify the relationship between these \emph{monoidal matrices} and $\mathcal{L}$-hypermatrices. We also show that, in the most favourable cases, the structure of the former can be simplified to a pair constituted by an algebra and a \emph{fuzzy subset} of its universe. We use these simplified matrices to provide a completeness theorem for a multiset-theoretic companion of infinite-valued \L{}ukasiewicz logic.

Finally, we introduce Hilbert systems suited for multiset consequence and prove that their derivability relations are substitution\=/invariant \mdr's; conversely, every substitution\=/invariant \mdr\ is shown to arise as the derivability relation of some such Hilbert system. As an example, we provide a Hilbert-style axiomatisation of the above-mentioned multiset-theoretic companion of infinite-valued \L{}ukasiewicz logic. The main result is:

\begin{theorem}[see Theorem~\ref{szko}]
Any substitution\=/invariant \mdr\ coincides with the derivability relation of some axiomatic system $\mathsf{AS}$.
\end{theorem}

\section{Deductive relations}\label{dudu}

\subsection{Basic definitions and facts}

We emphasised in our introduction that we need to consider more general relations than Galatos and Tsinakis' \gtcr's if we want to properly account for multiset consequence relations introduced in Definition~\ref{d:MDR}.
By Proposition~\ref{p:MutliPomon}, the set of finite multisets of formulas of a given language can be equipped with the structure of a dually integral Abelian pomonoid. This leads us to the next definition.

\begin{definition}\label{Def:DeductiveRelation}
A \emph{deductive relation} (\dr) on a dually integral Abelian pomonoid $\mathbf{R}=\langle R,\leq,+,0\rangle$ is a relation $\vdash$ on $R$ such that for every $a,b,c\in R$:
\begin{enumerate}
\item[$\bullet$] If\/ $a\leq b$, then $b\vdash a$. \hfill (\emph{Generalised Reflexivity})

\item[$\bullet$] If\/ $a\vdash b$ and $b\vdash c$, then $a\vdash c$. \hfill (\emph{Transitivity})

\item[$\bullet$] If\/ $a\vdash b$, then $a+c\vdash b+c$. \hfill (\emph{Compatibility})
\end{enumerate}

A \dr\ is \emph{finitary} if for each compact\footnote{An element $a\in R$ is \emph{compact} if for each directed set $D\subseteq R$ which has a supremum $\sup(D)  \geq a$ we have $d\geq a$ for some $d\in D$.} element $b$ such that $a\vdash b$ there is a compact element $a'\leq a$ such that $a'\vdash b$.
\end{definition}

\newpage

Observe that $\vdash$ is a compatible preordering of $R$ and $0$ is a $\vdash$\=/maximum. Also observe that, using properties of compatible preorderings on Abelian monoids and the dual integrality of $\mathbf{R}$, we have:

\begin{lemma}
Let\/ $\vdash$ be a \dr\ on $\mathbf{R}=\langle R,\leq,+,0\rangle$. For all $a,b,c,d\in R$:
\begin{enumerate}
\item $a\vdash b$ and $c+b\vdash d$ imply $c+a\vdash d$. \hfill (Cut)

\item $a\vdash b$ implies $c+a\vdash b$. \hfill (Monotonicity)
\end{enumerate}
\end{lemma}

Some examples of \dr's follow hereafter. Our prime motivating example, the \emph{multiset deductive relations},
will be thoroughly studied in Section~\ref{energade}, where appropriate
particular examples will be given. It is easy to observe that:

\begin{example}\label{tamburo}
Let $\mathcal{L}$ be a propositional language. Multiset deductive relations on $\mathcal{L}$ are exactly the \dr's on $\mathbf{Fm}_{\mathcal{L}}^{\flat}$.
\end{example}

Of course, standard \acr's (hence, in particular, \tcr's) give rise to instances of deductive relations:

\begin{example}\label{carriola}
Let $\vdash$ be an \acr\ on the set $A$ (see Definition~\ref{derkommissar}) and let
\[
 \mathbf{R}^{\wp(A)} = \langle \wp(A),\subseteq,\cup,\emptyset\rangle.
\]
Then $\vdash'$, where, for all $X,Y\subseteq A$, $X\vdash'Y$ iff $X\vdash a$ for all $a\in Y$, is a \dr\ on $\mathbf{R}^{\wp(A)}$.
\end{example}

In view of the previous example, the reader will be curious to figure out whether \dr's also generalise \gtcr's (see Definition~\ref{def:gtcr}). The answer, here, is less straightforward.

\begin{proposition}\label{GTCR-DR}
Let $\mathbf{L}=\langle L,\wedge,\vee\rangle $ be a complete lattice with induced order $\leq$ and bottom element $0$. Then the structure
\[
 \alg{L}' =\langle L,\leq ,\vee,0\rangle
\]
is a dually integral Abelian pomonoid and any (finitary) \gtcr\ on $\mathbf{L}$ is a (finitary) \dr\ on $\mathbf{L}'$.
\end{proposition}

\begin{proof}
The only non-trivial part is Compatibility. Assume that $a\vdash b$. Thus, $a\vee c\vdash b$ and $a\vee c\vdash c$ and so $\bigvee\{x : a\vee c\vdash x\}\vdash b\vee c$. By the defining condition on \gtcr's, it follows that $a\vee c\vdash b\vee c$. 
\end{proof}

The converse direction does not hold in general. The next proposition provides a class of explicit counterexamples, indeed, a very wide one, because any substitution\=/invariant \tcr\ $\vdash$ with a theorem (i.e., an element $a$ such $\emptyset\vdash a$) containing a variable has infinitely many theorems.

\begin{proposition}
Let $\vdash$ be a finitary \acr\ with infinitely many theorems\footnote{Actually, any \acr\ where there is a finite set with
infinitely many consequences would do the job.} on $A$ and let\/ $\vdash'$ be defined as in Example~\ref{carriola}.
Then the relation defined by:
\[
 X\Vdash Y\text{ if there is a finite } Y' \subseteq Y\text{ such that } Y\setminus Y'\subseteq X \text{ and } X\vdash' Y'
\]
is a \dr\ on $\mathbf{R}^{\wp(A)}$ but it is not a \gtcr\ on $\wp(A)$.
\end{proposition}

\begin{proof}
For our first claim, the only nontrivial condition to check is Transitivity. Assume that $X\Vdash Y$ and $Y\Vdash Z$, and let $Y'$ and $Z'$ be finite sets with the required properties. We claim that the finite set $\bar{Z}=Z'\cup(Y'\cap Z)$ witnesses $X\Vdash Z$. In fact:
\begin{enumerate}
\item[$\bullet$] Clearly $\bar Z$ is a finite subset of $Z$.

\item[$\bullet$] $X\vdash'\bar{Z}$, because ${\Vdash}\subseteq{\vdash'}$ implies that $X\vdash'Z$.

\item[$\bullet$] The final condition is obtained by the following chain:
\begin{align*}
Z\setminus(Z'\cup(Y'\cap Z))  &= (Z\setminus Z')\cap(Z\setminus(Y'\cap Z))=(Z\setminus Z')\cap(Z\setminus Y') \subseteq Y\cap(Z\setminus Y')\subseteq Y\setminus Y' \subseteq X.
\end{align*}
\end{enumerate}

To conclude the proof, it is easy to observe that for each theorem $a$ we have $\emptyset\Vdash a$, yet clearly $\emptyset\not \vdash \{a : \emptyset\Vdash a\}$.
\end{proof}

\newpage

The above result is not that surprising, as our definition of \dr\ only employs finitary operations, as opposed to the use of infinite suprema in \gtcr's. Our approach has the advantage that our background structure could be much smaller, as the following proposition shows.

\begin{proposition}
Let\/ $\mathbf{L}=\langle L,\wedge,\vee\rangle $ be a complete algebraic lattice with induced order $\leq$ and bottom element $0$. Then the structure
\[
 K(\mathbf{L}) =\langle K(L),\leq,\vee,0\rangle,
\]
where $K(L)  $ is the set of  compact elements of\/ $\mathbf{L}$, is a dually integral Abelian pomonoid. Moreover, there is bijective correspondence between \emph{finitary} \gtcr's on $\mathbf{L}$ and \dr's on $K(\mathbf{L})$.\footnote{Note that all \dr's on $K(\mathbf{L})$ are finitary.}
\end{proposition}

\begin{proof}
For a start, note that finitary \gtcr's on complete algebraic lattices are fully determined by their subrelations between compact elements. In fact, for $x\in L$, let
\[
 C_{x} = \{c\in K(L) : c\leq x\}.
\]

Thus $x=\bigvee C_{x}$. Now, take $x,y\in L$. We have that $x\vdash y$ iff for each $c\in C_{y}$ there is $x_{c}\in C_{x}$ such that $x_{c}\vdash c$. As in Proposition~\ref{GTCR-DR}, it is possible to show that the restriction of any \gtcr\ on $\mathbf{L}$ is a \dr\ on $K(\mathbf{L})$, and the previous observation entails that if two \gtcr's differ they also differ on compact elements.

Conversely, assume that $\vdash$ is a \dr\ and define:
\[
 x\vdash'y \quad\iff\quad
 \text{for each } c \in C_{y}\text{ there is }x_{c}\in C_{x}\text{ such that }x_{c}\vdash c.
\]
The only condition in need of a proof is that $x\vdash'p$, where $p = \bigvee\{a\in L : {x\vdash'a}\}$. Let $C=\{c\in K(L) : x\vdash'c\}$. If we show that $C=C_{p}$, we are done. For the nontrivial inclusion, assume that $c\in K(L)  $ and $c\leq p$. Then there are $c_{1},\dots,c_{n}$ such that $c\leq c_{1}\vee\dots\vee c_{n}$ and $x\vdash'c_{j}$ for all $j\leq n$. This means that for all $j\leq n$ there is $x_{j}\leq x$ such that $x_{j}\vdash c_{j}$. Using properties of \dr's we obtain $x_{1}\vee\dots\vee x_{n}\vdash c_{1}\vee\dots\vee c_{n}$ and so $x\vdash'c$, i.e., $c\in C$. This mapping is clearly one-one and an inverse to the previous one.
\end{proof}

The next example identifies a deductive relation on fuzzy sets. It is introduced to underscore the generality of our framework, but it will not be further discussed in the remainder of this paper.

\begin{example}
Let $\mathcal{L}$ be the language of infinite-valued \L{}ukasiewicz logic \L . The relation ${\vdash}\subseteq\lbrack0,1]^{\Fm_{\mathcal{L}}} \times \lbrack0,1]^{\Fm_{\mathcal{L}}}$ defined as:
\begin{align*}
\Gamma\vdash\Delta \quad\iff\quad
  &\text{for each } [0,1]\text{-valued evaluation $e$ we have: if } e(\psi)\geq\Gamma(\psi)\\
  &\text{for each }\psi\in\mathit{Fm}_{\mathcal{L}}, \text{then } e(\psi)\geq\Delta(\psi) \text{ for each } \psi\in\mathit{Fm}_{\mathcal{L}}
\end{align*}
is a \dr\ on
\[
 \mathbf{R}^{[0,1]^{\mathit{Fm}_{\mathcal{L}}}} =
 \langle [0,1]^{\mathit{Fm}_{\mathcal{L}}},\leq,\vee,\emptyset\rangle,
\]
where $\emptyset(\varphi)=0$ for all $\varphi\in\mathit{Fm}_{\mathcal{L}}$ and $\vee$ is pointwise supremum.
\end{example}

Among the basic notions of algebraic logic that need to be redefined in our new framework, one certainly finds the concepts of \emph{theory} and \emph{theorem}. Here, we must stray away to a certain extent from the received orthodoxy. In view of Example~\ref{carriola}, given a \dr\ $\vdash$ on $\mathbf{R} = \langle R,\leq,+,0\rangle$, one would expect a $\vdash$\=/theory to be an \emph{element} of $R$ with certain properties. In AAL, in fact, a theory is a deductively closed set of formulas. In particular, the theory generated by a set of formulas $X$ is the smallest deductively closed set of formulas that includes $X$ --- or else, the largest $Y$ such that $X\vdash Y$ --- and has the property that its subsets are exactly the consequences of $X$. However, in the case of \mdr's (Definition~\ref{d:MDR}), this would not work. Such a ``largest consequence'' need not always exist, because it could happen, for instance, that $\Gamma\vdash\Delta$ and $\Gamma\vdash\Pi$ but $\Gamma\not \vdash \Delta\vee\Pi$. Nevertheless, it makes sense to collect all the consequences of a given multiset $\Gamma$ of formulas into a set and view \emph{the set itself} as the deductive closure of $\Gamma$. Abstracting away from this particular example, we are led to the following definition (recall that every \dr\ is a preorder on $R$).

\begin{definition}\label{teoria}
Let\/ $\vdash$ be a \dr\ on $\mathbf{R}$. A\/ $\vdash$\=/\emph{theory} (or simply a \emph{theory}, when $\vdash$ is understood) is a $\vdash$\=/upset $T$ of $R$. By $\mathrm{Th}(\vdash)$ we denote the family of all\/ $\vdash$\=/theories.
\end{definition}

\begin{proposition}
Let\/ $\vdash$ be a \dr\ on $\mathbf{R}$. Then $\mathrm{Th}(\vdash)$ is a closure system on $R$ (namely, a family of subsets of $R$ that contains $R$ and is closed under arbitrary intersections).
\end{proposition}

\begin{proof}
Clearly, $R$ is a theory. Suppose that $\{T_{i}\}_{i\in I}$ is a nonempty family of theories, $a\in\bigcap\{T_{i}\}_{i\in I}$ and $a\vdash b$. Given an arbitrary $T_{j}$ ($j\in I$), $a\in T_{j}$, whence $b\in T_{j}$. It follows that $b\in\bigcap\{T_{i}\}_{i\in I}$.
\end{proof}

We denote by $\mathrm{Th}_{\vdash}(X)$ the smallest $\vdash$\=/theory containing $X$, and by $\mathrm{Th}^{p}(\vdash)$ the set of \emph{principal} theories of the form $\mathrm{Th}_{\vdash}(a)$, for $a\in R$, which is just the principal $\vdash$\=/upset generated by $a$. The subscripts in $\mathrm{Th}_{\vdash}$ will be dropped when the deductive relation is clear from the context.

Let us note the each $\vdash$\=/theory is a union of principal theories, a quite unusual feature from the point of view of the general theory of closure systems. Actually we can prove even more:

\begin{proposition}
Let\/ $\vdash$ be a \dr\ on $\mathbf{R}$. Then
\[
 \mathrm{Th}_{\vdash}(X)  =\bigcup\{\mathrm{Th}_{\vdash}(x) : x\in X\}.
\]
\end{proposition}

\begin{proof}
Clearly if $x\in X$ then $\mathrm{Th}_{\vdash}(x)\subseteq\mathrm{Th}_{\vdash}(X)$ and so one inclusion follows. To prove the converse one it suffices to show that $T = \bigcup\{\mathrm{Th}_{\vdash}(x) : x \in X\}$ is a $\vdash$\=/theory (since, as such, it will contain the smallest $\vdash$\=/theory including $X$). Assume that $a\in T$ and $a\vdash b$. Thus we know that $a\in\mathrm{Th}_{\vdash}(x)$ for some $x\in X$ and so Transitivity completes the proof.
\end{proof}

\begin{definition}
Let\/ $\vdash$ be a \dr\ on $\mathbf{R}=\langle R,\leq,+,0\rangle$. The element
$b\in R$ is a $\vdash$\=/\emph{theorem} if\/ $0\vdash b$.
\end{definition}

Observe that $0$ is always a theorem\footnote{There are \dr's for which $0$ is the \emph{only} theorem --- for example, any \dr\ that stems from a \emph{theoremless} \acr\ as in Example~\ref{carriola} has this property.} and the theorems are exactly the $\vdash$\=/maximal elements of $R$, since for any $a\in R$ we have $a\vdash 0$. Note that $\geq$, which can be seen as the least \dr\ on $\mathbf{R}$, has $0$ as the only theorem.

The reader will recall that one of the main advantages of the notion of \gtcr\ is the fact that the collection of theories of a given \gtcr\ $\vdash$ over a complete lattice is itself a complete lattice, which is furthermore determined by $\vdash$. We prove an analogous results for \dr's and \emph{principal} theories.

\begin{theorem}\label{DRandTheories}
Let $\vdash$ be a \dr\ on $\mathbf{R} = \langle R,\leq,+,0\rangle$. Let us define $+^{\vdash}$ on $\mathrm{Th}^{p}(\vdash)$ as
\[
 \mathrm{Th}_\vdash(x)+^{\vdash}\mathrm{Th}_\vdash(y) = \mathrm{Th}_\vdash(x+y).
\]
Then
\[
\mathbf{Th}_{\vdash} = \langle\mathrm{Th}^{p}(\vdash),\subseteq,+^{\vdash},\mathrm{Th}(0)\rangle
\]
is a dually integral pomonoid and the mapping $\mathrm{Th}_\vdash\colon R\to\mathrm{Th}^{p}(\vdash)$ is a surjective morphism.
\end{theorem}

\begin{proof}
The relation $\subseteq$ is clearly a partial order on $\mathrm{Th}^{p}(\vdash)$ with $\mathrm{Th}(0)$ as a bottom element. Moreover, $\mathrm{Th}$ is order preserving: in fact, by Generalised Reflexivity $a\in\mathrm{Th}(a)$ and so, by Transitivity, $a\vdash b$ iff $\mathrm{Th}(b)\subseteq \mathrm{Th}(a)$. Thus, if $a\leq b$ then, by Generalised Reflexivity, $\mathrm{Th}(a)\subseteq\mathrm{Th}(b)$.

We now sketch the proof of the fact that the operation $+^{\vdash}$ is well defined. To this end, consider $a,b,c,d\in R$ such that $\mathrm{Th}(a)=\mathrm{Th}(b)$ and $\mathrm{Th}(c)=\mathrm{Th}(d)$. In particular, $b\vdash a$ and $d\vdash c$ and so by Compatibility and Transitivity $d+b\vdash a+c$, whence $\mathrm{Th}(a+c)\subseteq\mathrm{Th}(b+d)$. The other inclusion is proved analogously. The fact that $\langle\mathrm{Th}^{p}(\vdash),+^{\vdash},\mathrm{Th}(0)\rangle$ is an Abelian monoid is obvious, so we only need to check that $+^{\vdash}$ is compatible with the order. Let $a,b,c\in R$ be such that $\mathrm{Th}(a)\subseteq\mathrm{Th}(b)$. Thus $b\vdash a$ and so $b+c\vdash a+c$. As a consequence,
\[
 \mathrm{Th}(a)+^{\vdash}\mathrm{Th}(c) = \mathrm{Th}(a+c) \subseteq
 \mathrm{Th}(b+c) = \mathrm{Th}(b)+^{\vdash}\mathrm{Th}(c).
\]
The fact that $\mathrm{Th}$ is a surjective morphism is again obvious.
\end{proof}

\subsection{Deductive operators and systems}

In AAL, propositional logics can be introduced in three different but equivalent ways: via consequence relations, via closure operators, and via closure systems. The same is true of the approach we have taken. While deductive relations are abstract counterparts of \tcr's, our next goal is to define suitable abstract notions of deductive operator and deductive system, in such a way as to generalise the lattice isomorphisms between the lattices of consequence relations, of closure operators, and of closure systems that are available in the traditional theory of AAL. Analogues of the classical concepts of closure operator and closure system can be defined as follows.

\begin{definition}\label{Def:DeductiveOper}
A \emph{deductive operator} (\doo) on a dually integral Abelian pomonoid $\mathbf{R} = \langle R,\leq,+,0\rangle$ is a map $\delta\colon R\to\mathcal{\wp}(R)$ such that for every $a,b,c\in R$:
\begin{enumerate}
\item[$\bullet$] $a\in\delta(a)$.\hfill (\emph{Enlargement})

\item[$\bullet$] If $a\leq b$, then $\delta(a)\subseteq\delta(b)$.\hfill (\emph{Order Preservation})

\item[$\bullet$] If $a\in\delta(b)$, then $\delta(a)\subseteq\delta(b)$.\hfill (\emph{Idempotency})

\item[$\bullet$] If $a\in\delta(b)$, then $a+c\in\delta(b+c)$.\hfill (\emph{Compatibility})
\end{enumerate}
\end{definition}

Observe that, in full analogy with Definition~\ref{teoria}, a \doo\ is a map from elements of $R$ to \emph{subsets} of $R$.

Next we define the notion of a deductive system;  recall that closure systems are systems of theories of some \acr, but as we have seen in the previous subsection the \emph{principal} theories are the crucial ones in our framework (as they can be seen as universes of dually integral Abelian pomonoids). This leads to the following definition:

\begin{definition}\label{Def:DeductiveSystem}
A \emph{deductive system} (\ds) on a dually integral Abelian pomonoid $\mathbf{R} = \langle R,\leq,+,0\rangle$ is a family $\mathcal{C}\subseteq\wp(R)$ of $\leq$-downsets of $R$ such that for the mapping $\delta_{\mathcal{C}}\colon R\to\wp(R)  $ defined by $\delta_{\mathcal{C}}(x)=\bigcap\{C\in\mathcal{C} : x\in C\}$ we have $\delta_{\mathcal{C}}(R)=\mathcal{C}$ and if $\delta_{\mathcal{C}}(x)\subseteq\delta_{\mathcal{C}}(y)$, then $\delta_{\mathcal{C}}(x+z)\subseteq\delta_{\mathcal{C}}(y+z)$ for all $x,y,z\in R$.
\end{definition}

The subscript, or superscript, $\mathcal{C}$ will be omitted whenever it is not needed to clarify potential confusions.

\begin{proposition}\label{prop:dsANDdo}
Let $\mathcal{C}$ be a \ds\ on a dually integral Abelian
pomonoid $\mathbf{R}$. Then the mapping $\delta_{\mathcal{C}}$ is a \doo\ on
$\mathbf{R}$.
\end{proposition}

\begin{proof}
Enlargement is obvious. For Idempotency, assume that $a\in\delta(b)$, $c\in\delta(a)$ and we have $X\in\mathcal{C}$ such that $b\in X$. Then $a\in X$ (due to the first assumption) and so $c\in X$ (due to the second assumption), i.e., $c\in\delta(b)$. If $a\leq b$, then (as each $X$ is $\leq$-downset) $a\in\delta(b)$ and so by Idempotency $\delta(a)\subseteq\delta(b)$, which takes care of Order Preservation. Finally, we prove Compatibility: using the conditions we already established, if $a\in\delta(b)$, then $\delta(a)\subseteq\delta(b)$ and so $\delta(a+c)\subseteq\delta(b+c)$, whence $a+c\in\delta(b+c)$.
\end{proof}

Given a dually integral Abelian pomonoid $\mathbf{R} = \langle R,\leq ,+,0\rangle$, we respectively denote by $\mathrm{\operatorname{Re} l}(\mathbf{R}),\mathrm{Op}(\mathbf{R})$ and $\mathrm{Sys}(\mathbf{R})$ the sets of deductive relations, deductive operators and deductive systems on $\mathbf{R}$. We next define partial orders on these sets. $\mathrm{\operatorname{Rel}}(\mathbf{R})$ will be viewed as partially ordered by set inclusion, and $\mathrm{Sys}(\mathbf{R})$ by supersethood. We define an order on $\mathrm{Op}(\mathbf{R})$ as follows: given $\delta,\gamma \in\mathrm{Op}(\mathbf{R})$, we set $\delta\preccurlyeq\gamma$ iff $\delta(a)\subseteq\gamma(a)$ for every $a\in R$.

\begin{theorem}\label{Thm:Trinity}
If\/ $\mathbf{R} = \langle R,\leq,+,0\rangle$ is a dually integral Abelian pomonoid, then the posets $\langle \mathrm{\operatorname{Rel}}(\mathbf{R}),\subseteq\rangle $, $\langle \mathrm{Op}(\mathbf{R}),\preccurlyeq\rangle $ and $\langle \mathrm{Sys}(\mathbf{R}),\supseteq\rangle $ are isomorphic.
\end{theorem}

\begin{proof}
We define maps $\delta_{()}\colon\mathrm{\operatorname{Re}l}(\mathbf{R})\to\mathrm{Op}(\mathbf{R})$ and $\vdash_{()}\colon\mathrm{Op}(\mathbf{R} \to\mathrm{\operatorname{Re}l}(\mathbf{R})$ as follows:
\[
 {\vdash_{\delta}} = \{\langle a,b\rangle : b\in\delta(a)\}
 \quad\text{and}\quad
 \delta_{\vdash}(a) = \mathrm{Th}_{\vdash}(a).
\]

It is easy to show that they are well defined and monotone. We now prove they are mutually inverse. Consider ${\vdash}\in\mathrm{Rel}(\mathbf{R})$ and observe that $x\vdash y$ iff $y\in\delta_{\vdash}(x)$ iff $x\vdash_{\delta_{\vdash}}y$. On the other hand, for any $\delta\in\mathrm{Op}(\mathbf{R})$:
\[
 x\in\delta(y) \quad\iff\quad
 y\vdash_{\delta}x \quad\iff\quad
 x\in\delta_{\vdash_{\delta}}(y).
\]
Next, we handle the slightly more complex case of \dr's and \ds's. We define two maps $\mathcal{C}_{()}\colon\mathrm{Rel}(\mathbf{R})\to\mathrm{Sys} (\mathbf{R})$ and $\vdash_{()}\colon\mathrm{Sys}(\mathbf{R})\to \mathrm{Rel}(\mathbf{R})$ as follows:
\[
 \mathcal{C}_{\vdash} = \{\mathrm{Th}_{\vdash}(a) : a\in R\}
 \quad\text{and}\quad
 {\vdash_{\mathcal{C}}} = \{\langle a,b\rangle : b\in\delta_{\mathcal{C}}(a)\}.
\]
Using Theorem~\ref{DRandTheories}, we conclude that $\mathcal{C}_{\vdash}$ is indeed a \ds, while Proposition~\ref{prop:dsANDdo} guarantees that $\delta_{\mathcal{C}}$ is a \doo. Monotonicity of both functions is obvious and so is the fact that they are mutually inverse. For the sake of completeness, we observe that $\mathcal{C}_{\delta}=\{\delta(a) : a\in R\}$ maps $\mathrm{Op}(\mathbf{R})$ to $\mathrm{Sys}(\mathbf{R})$, and that such a mapping is inverted by the mapping that sends any \ds\ $\mathcal{C}$ to the already defined deductive operator $\delta_{\mathcal{C}}$.
\end{proof}

\begin{corollary}\label{latop}
$\langle \mathrm{\operatorname{Rel}}(\mathbf{R}),\subseteq\rangle $, $\langle \mathrm{Op}(\mathbf{R}),\preccurlyeq\rangle $ and $\langle\mathrm{Sys}(\mathbf{R}),\supseteq\rangle$ are complete lattices.
\end{corollary}

\begin{proof}
By Theorem~\ref{Thm:Trinity}, it will suffice to prove our claim for any one of these posets, say $\langle \mathrm{Op}(\mathbf{R}),\preccurlyeq\rangle$. It is obvious that $\delta_{1}$, defined by
\[
 \delta_{1}(a)  =R\text{ for all }a\in R
\]
is a \doo\ and that it is the top element w.r.t.\ $\preccurlyeq$. In order to see that \doo's on $\mathbf{R}$ form a complete lattice, it is enough to see that for every family of \doo's $\{\delta_{i} : i\in I\}$, the map $\bigwedge_{i\in I}\delta_{i}$ defined by
\[
 \bigwedge_{i\in I}\delta_{i}(a) = \bigcap\{\delta_{i}(a) : i\in I\}
\]
is again a \doo. Enlargement is clear. For Order Preservation, suppose that $a\leq b$ and $c\in\delta_{i}(a)  $ for all $i\in I$. Then for all $i\in I$ we have that $\delta_{i}(a)\subseteq\delta_{i}(b)$, whence $c\in\bigcap\{\delta_{i}(b) : i\in I\}$. As to Idempotency, suppose $a\in\delta_{i}(b)  $ for all $i\in I$, and $c\in\delta _{i}(a)  $ for all $i\in I$. Then for all $i\in I$ we have that $\delta_{i}(a)  \subseteq\delta_{i}(b)$ and so again $c\in\bigcap\{\delta_{i}(b) : i\in I\}$. Compatibility, once more, is clear.
\end{proof}

\subsection{Blok--J\'{o}nsson companions of deductive relations}

Deductive relations, deductive operators and deductive systems respectively give rise to special kinds of \acr's (Definition~\ref{derkommissar}), closure operators and closure systems. In the present subsection, we point out the fact that there is a significant transfer of information from the original relations, operators and systems to these ``Blok--J\'{o}nsson companions'', which we now proceed to define.

Given a \dr\/ $\vdash$ on $\mathbf{R}=\langle R,\leq,+,0\rangle$, its \emph{Blok--J\'{o}nsson companion} is the relation ${\vdash^{BJ}}\subseteq \wp(R)  \times R$ defined as follows for every $X\subseteq R$ and every $a\in R$:
\[
 X\vdash^{\mathit{BJ}}a \quad\iff\quad \text{there is } y\in X \text{ s.t.\ } y\vdash a.
\]

\begin{lemma}\label{zenzero}
If\/ $\vdash$ is a \dr\ on $\mathbf{R}=\langle R,\leq,+,0\rangle$, then $\vdash^{\mathit{BJ}}$ is an \acr\ on $R$.
\end{lemma}

\begin{proof}
For Reflexivity, suppose $x\in X$; we want to show that $X$ $\vdash^{\mathit{BJ}}x$, i.e.\ that there is $y\in X$ s.t.\ $y\vdash x$. However, by the reflexivity of
$\vdash$, $x$ itself fits the bill.

Monotonicity is straightforward from the definition of $\vdash^{\mathit{BJ}}$.

For Cut, we want to show that $X$ $\vdash^{\mathit{BJ}}y$ and $Z\vdash^{\mathit{BJ}}x$ for all $x\in X$ imply that $Z\vdash^{\mathit{BJ}}y$. In fact, suppose there is $x\in X$ s.t.\ $x\vdash y$. This implies that there is $z\in Z$ s.t.\ $z\vdash x$, whence by transitivity of $\vdash$, $z\vdash y$, which in turns entails that $Z$ $\vdash^{\mathit{BJ}}y$.
\end{proof}

Along the same lines, deductive operators and deductive systems on $\mathbf{R}$ can be lifted to closure operators and closure systems, respectively, on the base set~$R$.

\begin{lemma}\label{clous}
Given a \doo\ $\delta$ on $\mathbf{R} = \langle R,\leq,+,0\rangle$, the map $\delta^{\mathit{BJ}}\colon \wp(R)\to\wp(R)$ defined as
\[
 \delta^{\mathit{BJ}}(X)  = \bigcup \{\delta(x) : x\in X\}
\]
is a closure operator.
\end{lemma}

\begin{proof}
First, if $a\in X\subseteq R$, then $a\in\delta(a)\subseteq\delta^{\mathit{BJ}}(X)$, and therefore $X\subseteq\delta^{\mathit{BJ}}(X)$. If $X\subseteq Y$, then obviously $\delta^{\mathit{BJ}}(X)\subseteq\delta^{\mathit{BJ}}(Y)$. And finally, assume that $z\in \delta^{\mathit{BJ}}(\delta^{\mathit{BJ}}(X))  $, i.e., there are $y,x$ such that $x\in X,y\in\delta(x),z\in\delta(y)$. By Idempotency and Generalised Reflexivity, $z\in\delta(z)\subseteq\delta(y)\subseteq\delta(x)$, which means $z\in\delta^{\mathit{BJ}}(X)$.
\end{proof}

\begin{lemma}\label{gummo}
Given a \ds\/ $\mathcal{C}$ on $\mathbf{R} = \langle R,\leq,+,0\rangle $, the family $\mathcal{C}^{\mathit{BJ}}=\{\bigcup \mathcal{Y} : \mathcal{Y}\subseteq\mathcal{C}\}$ is a closure system.
\end{lemma}

\begin{proof}
Recall that $\mathcal{C}=\{\delta_{\mathcal{C}}(x) : x\in R\}$ and as $x\in\delta_{\mathcal{C}}(x)$ then $\bigcup\mathcal{C}=R$. Now we only have to prove that for $\mathcal{X}\subseteq\mathcal{C}$ we have $\bigcap \mathcal{X}\in\mathcal{C}$. If $x\in\bigcap\mathcal{X}$ then for each $X\in\mathcal{X}$ there is $c_{X}$ such that $x\in\delta_{\mathcal{C}} (c_{X})\subseteq X$, thus also $\delta_{\mathcal{C}}(x)\subseteq \delta_{\mathcal{C}}(c_{X})\subseteq X$ and so $\delta_{\mathcal{C} }(x)\subseteq\bigcap\mathcal{X}$. To conclude the proof just observe that
\[
 \bigcap\mathcal{X} =
 \bigcup\Bigl\{\delta_{\mathcal{C}}(x) : x\in\bigcap\mathcal{X}\Bigr\} \in \mathcal{C}^{\mathit{BJ}}.
\]
\end{proof}

In full analogy with the above, we call $\delta^{\mathit{BJ}}$ and $\mathcal{C}^{\mathit{BJ}}$ the \emph{Blok--J\'{o}nsson companions} of $\delta$ and $\mathcal{C}$, respectively. Observe that:

\begin{lemma}\label{gummi}
Let\/ $\vdash$ be a \dr\ on $\mathbf{R} = \langle R,\leq ,+,0\rangle $. The $\vdash$\=/theories are the theories (in the sense of Blok--J\'{o}nsson) of\/ $\vdash^{\mathit{BJ}}$. Namely, for $T\subseteq R$, t.f.a.e.:
\begin{enumerate}
\item $T$ is a $\vdash$\=/upset of $R$.

\item $T\vdash^{\mathit{BJ}}x$ implies $x\in T$.
\end{enumerate}
\end{lemma}

\begin{proof}
Let us first suppose $T$ is a $\vdash$\=/upset of $R$ and there is $y\in T$ such that\ $y\vdash x$. Then $x\in T$. Conversely, suppose that for any $x$, $T$ $\vdash^{\mathit{BJ}}x$ implies $x\in T$, that $y\in T$, and that $y\vdash z$. So $T$ $\vdash^{\mathit{BJ}}z$, whence $z\in T$.
\end{proof}

Let us continue to use the notations $\vdash_{\delta},\vdash_{\mathcal{C}},\mathcal{C}_{\vdash},\mathcal{C}_{\delta},\delta_{\vdash}$, and $\delta_{\mathcal{C}}$ for the correspondences between \dr's, \doo's and \ds's on $\mathbf{R}$ spelt out in Theorem~\ref{Thm:Trinity}.
With an innocent notational abuse, we employ the same symbols for the standard correspondences between the sets $\mathrm{Acr}(R)$ of \acr's, closure operators $\mathrm{Clop}(R)$ and closure systems $\mathrm{Clos}(R)$, all on $R$. We now prove that the relation of ``taking the Blok--J\'{o}nsson companion'' commutes with these functions.

\begin{theorem}\label{gorgonzola}
The following diagrams
\[
\begin{tikzcd}
\mathrm{Rel}(\mathbf{R}) \ar[r, "\delta_{\vdash}"]\ar[d, swap, "\vdash^{\mathit{BJ}}"]
  & \mathrm{Op}(\mathbf{R})\ar[d, "\delta^{\mathit{BJ}}"]\\
\mathrm{Acr}(R) \ar[r, "\delta_{\vdash}"] & \mathrm{Clop}(R)
\end{tikzcd}
\qquad
\begin{tikzcd}
\mathrm{Rel}(\mathbf{R}) \ar[r, "\mathcal{C}_{\vdash}"]\ar[d, swap, "\vdash^{\mathit{BJ}}"]
   & \mathrm{Sys}(\mathbf{R})\ar[d, "\mathcal{C}^{\mathit{BJ}}"]\\
 \mathrm{Acr}(R)\ar[r, "\mathcal{C}_{\vdash}"] & \mathrm{Clos}(R)
\end{tikzcd}\]
\[
\begin{tikzcd}
\mathrm{Op}(\mathbf{R})\ar[r, "\delta_{\mathcal{C}}"]\ar[d, swap, "\delta^{\mathit{BJ}}"]
   & \mathrm{Sys}(\mathbf{R})\ar[d, "\mathcal{C}^{\mathit{BJ}}"]\\
\mathrm{Clop}(R)\ar[r,"\delta_{\mathcal{C}}"] & \mathrm{Clos}(R)
\end{tikzcd}
\]
as well as the ones we obtain by reversing the above correspondences, are all commutative.
\end{theorem}

\begin{proof}
The correspondences are well defined by Theorem~\ref{Thm:Trinity} and Lemmas~\ref{zenzero}, \ref{clous}, and~\ref{gummo}.
We now take care of some of the commutations. We show that ${(\vdash_{\delta}) ^{\mathit{BJ}}}={\vdash_{\delta^{\mathit{BJ}}}}$; the other commutations are established similarly. In fact, $\langle X,a\rangle \in(\vdash_{\delta})  ^{\mathit{BJ}}$ iff there exists $x\in X$ s.t.\ $x\vdash_{\delta}a$, which in turn holds iff there exists $x\in X$ s.t.\ $a\in\delta(x)  $. But this just means that $a\in \bigcup\nolimits_{x\in X}\delta(x)  $, which amounts to $X\vdash_{\delta^{\mathit{BJ}}}a$. Similarly, $(\delta_{\vdash}) ^{\mathit{BJ}}=\delta_{\vdash^{\mathit{BJ}}}$. In fact,
\begin{flalign*}
&&
  (\delta_{\vdash})  ^{\mathit{BJ}}(X)   &= \bigcup\nolimits_{a\in X}\mathrm{Th}_{\vdash}(a) 
  = \{b\in R : a\vdash b\text{ for some }a\in X\}\\
&&
  &= \{b\in R : X\vdash^{\mathit{BJ}}b\} 
   = \delta_{\vdash^{\mathit{BJ}}}(X).
\end{flalign*}
\end{proof}

This theorem implies, in particular, the following corollary:

\begin{corollary}\label{fortezza}
Let $\mathbf{R}=\langle R,\leq,+,0\rangle $ be a
dually integral Abelian pomonoid. The complete lattices of Blok--J\'{o}nsson companions of \dr's, \doo's, and \ds's on $\mathbf{R}$ are isomorphic.
\end{corollary}

\section{Action-invariance}\label{azione}

One of the remarkable achievements of Blok and J\'{o}nsson's treatment of logical consequence is its purely abstract account of substitution\=/invariance. Resorting to appropriate monoidal actions, Blok and J\'{o}nsson effectively sidestep the problem brought about by their use of sets with no structure whatsoever to be preserved. As we have seen, Galatos and Tsinakis turn this insight into the starting point for their categorical foundation of the whole subject. It would be highly desirable, then, to lay down a comparable treatment of action-invariance in our framework. This will be done by equipping our Abelian pomonoids with appropriate monoidal actions.

Our guiding example will again be given by multiset deductive relations, i.e., \dr's on $\mathbf{Fm}_{\mathcal{L}}^{\flat}$, as we can naturally call an \mdr\ $\vdash$ \emph{substitution\=/invariant} if for every $\mathcal{L}$-substitution $\sigma$ and for every $\Gamma,\Delta\in\mathit{Fm}_{\mathcal{L}}^{\flat}$:
\[
 \Gamma\vdash\Delta \quad\then\quad \sigma(\Gamma) \vdash\sigma(\Delta) .
\]

\subsection{A categorical setting}

For a start, let us recall the notion of partially ordered semiring~\cite[Ch.~3]{HW}.

\begin{definition}
A \emph{partially ordered semiring}, or \emph{po\=/semiring}, is a structure $\mathbf{A} = \langle A,\leq,+,\cdot,0,1\rangle$ where:
\begin{enumerate}
\item $\langle A,\cdot,1\rangle$ is a monoid.

\item $\langle A,\leq,+,0\rangle$ is an Abelian pomonoid.

\item $\sigma\cdot0=0\cdot\sigma=0$ for all $\sigma\in A$;

\item For every $\sigma,\pi,\varepsilon\in A$, we have
\[
 \pi\cdot(\sigma+\varepsilon) = (\pi\cdot\sigma)+(\pi\cdot\varepsilon)
 \quad\text{and}\quad
 (\sigma+\varepsilon)\cdot\pi = (\sigma\cdot\pi)+(\varepsilon\cdot\pi).
\]

\item If $\sigma\leq\pi$ and $0\leq\varepsilon$, then $\sigma\cdot \varepsilon\leq\pi\cdot\varepsilon$ and $\varepsilon\cdot\sigma\leq \varepsilon\cdot\pi$.
\end{enumerate}
\end{definition}

A po\=/semiring $\mathbf{A}=\langle A,\leq,+,\cdot,0,1\rangle$ is \emph{dually integral} iff $\langle A,\leq,+,0\rangle$ is dually integral as a pomonoid. Of course, the dual integrality condition ``kills'' many among the interesting examples of po\=/semirings, including all nontrivial po\=/rings.

Our chief example of dually integral po\=/semiring will be the semiring of finite multisets of substitutions on formulas of a propositional language $\mathcal{L}$. The role it will play here is analogous to the role played in Galatos and Tsinakis' theory by the complete residuated lattice of sets of $\mathcal{L}$-substitutions.

\begin{example}\label{Exa:Sigma}
Let $\mathcal{L}$ be a propositional language, and let $\textup{End}(\mathbf{Fm}_{\mathcal{L}})$ be the set of  \emph{substitutions} of\/ $\mathbf{Fm}_{\mathcal{L}}$. The structure
\[
 \boldsymbol{\Sigma}_{\mathcal{L}} =
 \langle\textup{End}(\mathbf{Fm}_{\mathcal{L}})^{\flat},\leqslant,\uplus,\cdot,\emptyset,[id_{\Fm_{\mathcal{L}}}]\rangle,
\]
where, for $\mathfrak{X}=[  \sigma_{1},\dots,\sigma_{n}]$, $\mathfrak{Y}=[  \pi_{1},\dots,\pi_{m}]  \in \textup{End}(\mathbf{Fm}_{\mathcal{L}})^\flat$,
\[
 \mathfrak{X}\cdot\mathfrak{Y} =
 [\sigma_{1}\circ\pi_{1},\dots,\sigma_{1}\circ\pi_{m},\dots,\sigma_{n}\circ\pi_{1},\dots,\sigma_{n}\circ\pi_{m}],
\]
is a dually integral po\=/semiring.
\end{example}

With this notion in our quiver, in order to get going we only need to endow our dually integral Abelian pomonoids from the previous section with a suitable operation of multiplication by a scalar.

\begin{definition}
Let $\mathbf{A} = \langle A,\leq^{\mathbf{A}},+^{\mathbf{A}},\cdot ^{\mathbf{A}},0^{\mathbf{A}},1^{\mathbf{A}}\rangle $ be a dually integral po\=/semiring. An $\mathbf{A}$-\emph{module} is a structure $\boldsymbol{R} = \langle R,\leq^{\boldsymbol{R}}, +^{\boldsymbol{R}}, 0^{\boldsymbol{R}},\ast^{\boldsymbol{R}}\rangle $ where $\langle R,{\leq^{\boldsymbol{R}},} +^{\boldsymbol{R}},0^{\boldsymbol{R}}\rangle $ is a dually integral Abelian pomonoid and $\ast^{\boldsymbol{R}}\colon A\times R\to R$ is an action of $\langle A,\cdot^{\mathbf{A}},1^{\mathbf{A}}\rangle $ on $R$ that is order-preserving in both coordinates and distributes over $+^{\boldsymbol{R}}$. In symbols:
\begin{enumerate}
\item $(\sigma\cdot^{\mathbf{A}}\pi)\ast^{\boldsymbol{R}}a = \sigma\ast^{\boldsymbol{R}}(\pi\ast^{\boldsymbol{R}}a)$

\item $1^{\boldsymbol{A}}\ast^{\boldsymbol{R}}a=a$

\item $0^{\mathbf{A}}\ast^{\boldsymbol{R}}a=0^{\mathbf{R}}$

\item $(\sigma\ast^{\boldsymbol{R}}a)+^{\boldsymbol{R}}(\sigma\ast^{\boldsymbol{R}}b) = \sigma\ast^{\boldsymbol{R}}(a+^{\boldsymbol{R}}b)$

\item $(\sigma+^{\mathbf{A}}\pi)\ast^{\boldsymbol{R}}a = (\sigma\ast^{\boldsymbol{R}}a)+^{\boldsymbol{R}}(\pi\ast^{\boldsymbol{R}}a)$

\item If\/ $\sigma\leq^{\mathbf{A}}\pi$, then $\sigma\ast^{\boldsymbol{R}} a\leq^{\boldsymbol{R}}\pi\ast^{\boldsymbol{R}}a$

\item If\/ $a\leq^{\mathbf{R}}b$, then $\sigma\ast^{\boldsymbol{R}} a\leq^{\boldsymbol{R}}\sigma\ast^{\boldsymbol{R}}b$.
\end{enumerate}
\end{definition}

\begin{example}\label{Exa:FinMultisets}
Consider the po\=/semiring $\boldsymbol{\Sigma}_{\mathcal{L}}$ defined in Example~\ref{Exa:Sigma}, and let $\Mult_{\mathcal{L}}$ $=\langle \Fm_{\mathcal{L}}^{\flat},\leqslant,\uplus,\emptyset,\ast\rangle $, where for
\[
 \mathfrak{X} = [\sigma_{1},\dots,\sigma_{n}] \in \textup{End}(\mathbf{Fm}_{\mathcal{L}})^\flat
 \quad\text{and}\quad
 \Gamma \in \Fm_{\mathcal{L}}^{\flat} ,
\]
we set, resorting to our usual notational conventions,
\[
\mathfrak{X}\ast\Gamma= \sigma_{1}(\Gamma) \uplus\dots\uplus\sigma_{n}(\Gamma).
\]
Then $\Mult_{\mathcal{L}}$ is a $\boldsymbol{\Sigma}_{\mathcal{L}}$\=/module.
\end{example}

Modules over a dually integral po\=/semiring can be naturally equipped with arrows as follows:

\begin{definition}
Let $\mathbf{A}$ be a dually integral po\=/semiring, and $\boldsymbol{R}$ and $\boldsymbol{S}$ be a pair of $\mathbf{A}$\=/modules. A \emph{morphism} $\tau\colon\boldsymbol{R}\to\boldsymbol{S}$ is a pomonoid homomorphism (i.e., an order\=/preserving monoid homomorphism) such that $\tau(\sigma\ast^{\boldsymbol{R}}a)=\sigma \ast^{\boldsymbol{S}}\tau(a)$ for every $\sigma\in A$ and $a\in R$.
\end{definition}

Given a dually integral po\=/semiring $\mathbf{A}$, the collection of $\mathbf{A}$\=/modules with morphisms between them forms a category in which composition and identity arrows are, respectively, standard composition of functions and identity functions. We denote this category by $\Mod{\mathbf{A}}$. Isomorphisms in the category $\Mod{\mathbf{A}}$ are precisely bijective morphisms that reflect the order.

From now on we will assume that $\mathbf{A}$ is a fixed, but otherwise arbitrary, dually integral po\=/semiring.

\begin{example}
It is expedient to remark that the setting of modules over complete residuated lattices can subsumed under the present one as follows. Recall that every complete residuated lattice $\mathbf{M}=\langle M,\wedge,\vee ,\cdot,\backslash,/,1\rangle $ can be naturally turned into a dually integral po\=/semiring $\mathcal{U}(\mathbf{M})=\langle M,\leq,\vee ,\cdot,0,1\rangle $ where $\leq$ and $0$ are respectively the order and the bottom element of the lattice reduct of $\mathbf{M}$. Then observe that every $\mathbf{M}$\=/module $\boldsymbol{L}=\langle L,\wedge,\vee,\star\rangle$ gives rise to a $\mathcal{U}(\mathbf{M})$\=/module $\mathcal{U}(\boldsymbol{L}) = \langle L,\leq,\vee,0,\star\rangle$, where $\leq$ and $0$ are respectively the order and the bottom element of the lattice reduct of $\boldsymbol{L}$. Finally, every translator $f\colon\boldsymbol{L}_{1}\to\boldsymbol{L}_{2}$ between $\mathbf{M}$\=/modules induces a morphism $\mathcal{U}(f) \colon\mathcal{U}(\boldsymbol{L}_{1})\to\mathcal{U}(\boldsymbol{L}_{2})$ of $\mathcal{U}(\mathbf{M})$\=/modules by setting $\mathcal{U} (f)(a)=f(a)$ for every $a\in L_{1}$. Summing up, the application $\mathcal{U}(\cdot)$ can be regarded as a forgetful functor from $\Mod{\mathbf{M}}$ to $\Mod{\mathcal{U}(\mathbf{M})}$, which reduces modules over a complete residuated lattice to modules over a dually integral po\=/semiring.
\end{example}

We are now ready to give an abstract formulation of action-invariant \dr's. Against the backdrop of Theorem~\ref{Thm:Trinity}, these deductive relations can be presented equivalently as deductive operators or as deductive systems. As a matter of fact, it turns out that working with \doo's is more convenient, although similar definitions and results can be obtained by putting the other two concepts to the forefront.

\begin{definition}\label{Def:AbstractSenseOfStructural}
An \emph{action-invariant} \doo\ on an $\mathbf{A}$\=/module $\boldsymbol{R}$ is a \doo\ $\delta$ on its pomonoid reduct $\langle R,\leq,+,0\rangle$ such that for every $\sigma\in A$ and $a,b\in R$:
\[
 a\in\delta(b) \quad\then\quad \sigma\ast a\in\delta(\sigma\ast b).
\]
\end{definition}

To exemplify this concept, we point out that substitution\=/invariant \mdr's give rise to deductive operators that are action-invariant according to the definition just given.

\begin{proposition}
Let $\mathcal{L}$ be a propositional language. Then an \mdr\ $\vdash$ on $\mathcal{L}$ is substitution\=/invariant iff\/ $\delta_{\vdash}$ is an action-invariant \doo\ on the $\boldsymbol{\Sigma}_{\mathcal{L}}$\=/module $\Mult_{\mathcal{L}}$. Similarly, $\delta$ is an action-invariant \doo\ on $\Mult_{\mathcal{L}}$ iff\/ $\vdash_{\delta}$ is a substitution\=/invariant \mdr\ on $\mathcal{L}$.
\end{proposition}

\begin{proof}
First consider an \mdr\ $\vdash$ on $\mathcal{L}$. Suppose that $\vdash$ is substitution\=/invariant. As $\mathbf{Fm}_{\mathcal{L}}^{\flat}$ is a reduct of the $\boldsymbol{\Sigma}_{\mathcal{L}}$\=/module $\Mult_{\mathcal{L}}$, by Theorem~\ref{Thm:Trinity} we only have to show that $\delta_{\vdash}$ is action-invariant. Assume that $\Delta\in\delta_{\vdash}(\Gamma)$ and consider $\mathfrak{X} = [\sigma_{1},\dots,\sigma_{k}] \in\textup{End}(\mathbf{Fm}_{\mathcal{L}})^{\flat}$. From the definition of $\delta_{\vdash}$ and substitution\=/invariance of $\vdash$ it follows that, for every $i\leq k$, we have that $\sigma_{i}(\Gamma)  \vdash \sigma_{i}(\Delta)  $. Now, applying Compatibility several times, we obtain that
\[
 \biguplus_{i\leq k}  \sigma_{i}(\Gamma)  \vdash\biguplus_{i\leq k} \sigma_{i}(\Delta).
\]
The above display amounts exactly to the fact that $\mathfrak{X}\ast\Delta\in \delta_{\vdash}(\mathfrak{X}\ast\Gamma)$. Hence, we conclude that $\delta_{\vdash}$ is action-invariant according to Definition~\ref{Def:AbstractSenseOfStructural}.

Conversely, suppose that $\delta_{\vdash}$ is action-invariant and that $\Gamma\vdash\Delta$. Consider a substitution $\sigma$. First observe that $\Delta\in\delta(\Gamma)$. Since $\delta_{\vdash}$ is action-invariant and $[  \sigma]  \in\textup{End}(\mathbf{Fm}_{\mathcal{L}})^{\flat}$, we have that $[\sigma]  \ast\Delta\in\delta_{\vdash}(\sigma \ast\Gamma)$, which means exactly $\sigma(\Gamma)  \vdash \sigma(\Delta)$. Hence we conclude that $\vdash$ is substitution\=/invariant.

The second claim follows from the first one, together with $\delta_{\vdash_{\delta}}=\delta$ (Theorem~\ref{Thm:Trinity}).
\end{proof}

Given an action-invariant \doo\ $\delta$ on $\boldsymbol{R}$, we define a structure
\[
 \boldsymbol{R}_{\delta} = \langle\delta(R),\subseteq,+^{\delta},\delta(0),\ast^{\delta}\rangle,
\]
where for every $\sigma\in A$ and $a,b\in R$:  $\delta(a)+^{\delta}\delta(b) = \delta(a+b)$ and  $\sigma\ast^{\delta}\delta(a) = \delta(\sigma\ast a)$.

\begin{lemma}
Let $\delta$ be an action-invariant \doo\ on the $\mathbf{A}$\=/module $\boldsymbol{R}$. Then $\boldsymbol{R}_{\delta}$ is a well\=/defined $\mathbf{A}$\=/module and the map $\delta\colon\boldsymbol{R}\to\boldsymbol{R}_{\delta}$ is a morphism.
\end{lemma}

\begin{proof}
Using Theorems~\ref{Thm:Trinity} and~\ref{DRandTheories}, we know that $\boldsymbol{R}_{\delta} = \langle\delta(R),\subseteq,+^{\delta},\delta(0)\rangle$
is a well\=/defined dually integral pomonoid. Now we show that the action $\ast^{\delta}$ is well defined too. Consider $\sigma\in A$ and $a,b\in R$ such that $\delta(a)=\delta(b)$. In particular, we have that $a\in\delta(b)$. By the action-invariance of $\delta$, we obtain that $\sigma\ast a\in \delta(\sigma\ast b)$. Thus we conclude that $\delta(\sigma\ast a)\subseteq \delta(\sigma\ast b)$. The other inclusion is proved analogously.

Next, we turn to prove that $\boldsymbol{R}_{\delta}$ is an $\mathbf{A}$\=/module. It only remains to establish the conditions regarding the action
$\ast^{\delta}$. It is clear that $\ast^{\delta}$ is order-preserving on the
first coordinate. We prove that the same holds for the second one. Consider
$\sigma\in A$ and $a,b\in R$ such that $\delta(a)\leq^{\delta}\delta(b)$. From
the action-invariance of $\delta$ it follows that $\sigma\ast a\in
\delta(\sigma\ast b)$ and, therefore, that $\delta(\sigma\ast a)\leq^{\delta}\delta(\sigma\ast b)$. We conclude that
\[
 \sigma\ast^{\delta}\delta(a) = \delta(\sigma\ast a)\leq^{\delta}\delta(\sigma\ast b)
 = \sigma\ast^{\delta}\delta(b).
\]

The fact that $\ast^{\delta}$ is a monoidal action and the distributivity conditions are easy exercises.
Finally, we prove that the map $\delta\colon\boldsymbol{R}\to \boldsymbol{R}_{\delta}$ is a morphism. Due to Theorem~\ref{DRandTheories}, it remains to show that $\delta$ respects the monoidal action, which follows directly from the definition of~$\ast^{\delta}$.
\end{proof}

We conclude this subsection be defining two maps which will play an important role in the next subsection.

\begin{lemma}\label{Lem:Technical}
Let $f \colon \boldsymbol{R} \to \boldsymbol{S}$ be a morphism between $\mathbf{A}$\=/modules.
\begin{enumerate}
\item  The map $f^{\ast}\colon R\to\mathcal{\wp}(R)$ defined as:
\[
 f^{\ast}(a) = f^{-1}(\{x : x \leq f(a)\})
\]
is an action-invariant \doo\ on $\boldsymbol{R}$.

\item  The map $\hat{f}\colon\boldsymbol{R}_{f^{\ast}}\to f[\boldsymbol{R}]$ defined as:
\[
 \hat{f}(f^{\ast}(a)) = f(a)
\]
is a well\=/defined isomorphism.
\end{enumerate}
\end{lemma}

\begin{proof} To prove the first claim, the only condition in Definition~\ref{Def:DeductiveOper} that stands in need of a check is Compatibility. Let $a,b,c\in R$ and suppose that $a\in f^{\ast}(b)$. This means that $f(a)\leq f(b)$. In particular, we have that
\[
 f(a+^{\boldsymbol{R}}c)
 = f(a)+^{\boldsymbol{S}}f(c)\leq^{\boldsymbol{S}}f(b)+^{\boldsymbol{S}}f(c)
 = f(b+^{\boldsymbol{R}}c).
\]
Hence we conclude that $a+c\in f^{\ast}(b+c)$. This shows that $f^{\ast}$ is a \doo. It remains to be shown that it is action-invariant. Consider $\sigma\in A$ and suppose that $a\in f^{\ast}(b)$. Then
\[
 f(\sigma\ast^{\boldsymbol{R}}a)
 = \sigma\ast^{\boldsymbol{S}}f(a)\leq^{\boldsymbol{S}}\sigma\ast^{\boldsymbol{S}}f(b)
 = f(\sigma\ast^{\boldsymbol{R}}b).
\]
Thus $\sigma\ast^{\boldsymbol{R}}a\in f^{\ast}(\sigma\ast^{\boldsymbol{R}}b)$, whence our conclusion follows.

To prove the second claim observe that the map $\hat{f}$ is well defined, since $\leq^{\boldsymbol{S}}$ is antisymmetric. It is clear that $\hat{f}$ is a bijection. Since isomorphisms in $\Mod{\mathbf{A}}$ are bijective morphisms, it suffices to prove that $\hat{f}$ is a morphism. But this is an exercise, using the definition of $\boldsymbol{R}_{f}$ and the fact that $f$ is a morphism.
\end{proof}

\subsection{Action-invariant representations}

The main result in~\cite{GT}, reproduced above as Theorem~\ref{calabro}, is an elegant and purely categorical characterisation of the modules over a complete residuated lattice for which an analogue of the Syntactic Isomorphism Theorem (Theorem~\ref{Thm:SyntacticIsoThm}) for algebraisable logics holds. The aim of this subsection is to obtain a similar result in the setting of modules over a dually integral po\=/semiring.

\begin{definition}
Let $\delta$ and $\gamma$ be action-invariant \doo's on the $\mathbf{A}$\=/modules $\boldsymbol{R}$ and $\boldsymbol{S}$, respectively.
\begin{enumerate}
\item An \emph{action-invariant representation} of $\delta$ into $\gamma$ is a morphism $\Phi\colon\boldsymbol{R}_{\delta}\to\boldsymbol{S}_{\gamma}$ that is injective and reflects the order.

\item A representation $\Phi$ of $\delta$ into $\gamma$ is \emph{induced} if there is a morphism $\tau\colon\boldsymbol{R}\to\boldsymbol{S}$ that makes the following diagram commute:
\[
\begin{tikzcd}
\boldsymbol{R} \ar[r, "\tau"] \ar[d,swap,"\delta"] & \boldsymbol{S} \ar[d, "\gamma"]\\
 \boldsymbol{R}_{\delta}\ar[r, "\Phi"] & \boldsymbol{S}_{\gamma}
\end{tikzcd}
\]

\item $\delta$ and $\gamma$ are \emph{equivalent} if the $\mathbf{A}$\=/modules $\boldsymbol{R}_{\delta}$ and $\boldsymbol{S}_{\gamma}$ are isomorphic.
\end{enumerate}
\end{definition}

\begin{definition}
An $\mathbf{A}$\=/module $\boldsymbol{R}$ has the  \emph{representation property} (REP) if for any other $\mathbf{A}$\=/module $\boldsymbol{S}$ and action-invariant \doo's $\delta$ and $\gamma$ on $\boldsymbol{R}$ and $\boldsymbol{S}$ respectively, every action-invariant representation of $\delta$ into $\gamma$ is induced.
\end{definition}

We are now ready to provide a characterisation of $\mathbf{A}$\=/modules with the REP in the spirit of Theorem~\ref{calabro}.

\begin{theorem}\label{nacchera}
An $\mathbf{A}$\=/module has the REP iff it is onto\=/projective
in $\Mod{\mathbf{A}}$.
\end{theorem}

\begin{proof}
The backbone of our argument is essentially the same as in~\cite[Lemma~8.1]{GT}. It is clear that every projective $\mathbf{A}$\=/module has the REP. Now, let $\boldsymbol{R}$ be an $\mathbf{A}$\=/module with the REP, and consider two morphisms $f\colon\boldsymbol{S}\to\boldsymbol{T}$ and $g\colon\boldsymbol{R}\to\boldsymbol{T}$ with $f$ onto. By Lemma~\ref{Lem:Technical}, the derived maps $f^{\ast}$ and $g^{\ast}$ are action-invariant \doo's on $\boldsymbol{R}$ and $\boldsymbol{S}$, respectively. Observe that $g(R) $ is the universe of a submodule $g(\boldsymbol{R})$ of $\boldsymbol{T}$. By Lemma~\ref{Lem:Technical} the following maps are isomorphisms:
\[
 \widehat{f}\colon\boldsymbol{S}_{f^{\ast}}\to\boldsymbol{T}
 \quad\text{and}\quad
 \widehat{g}\colon\boldsymbol{R}_{g^{\ast}}\to g(\boldsymbol{R}).
\]
Let $i\colon g(\boldsymbol{R})\to\boldsymbol{T}$ be the morphism given by the inclusion relation. Clearly the composition
\[
\hat{f}^{-1}\circ i\circ\hat{g}\colon\boldsymbol{R}_{g^{\ast}}\to \boldsymbol{S}_{f^{\ast}}
\]
is a representation of $g^{\ast}$ into $f^{\ast}$. Thus we can apply the fact that $\boldsymbol{R}$ has the REP, obtaining a morphism $h\colon\boldsymbol{R}\to\boldsymbol{S}$ such that
\[
 \hat{f}^{-1}\circ i\circ\hat{g}\circ g^{\ast} = f^{\ast}\circ h.
\]
Hence for every $a\in R$, we have that
\begin{align*}
f\circ h(a)  &= (\hat{f}\circ f^{\ast})\circ h(a)=\hat{f}\circ(f^{\ast}\circ h)(a)
  = \hat{f}\circ(\hat{f}^{-1}\circ i\circ\hat{g}\circ g^{\ast})(a)\\
  &= (\hat{f}\circ\hat{f}^{-1})\circ i\circ\hat{g}\circ g^{\ast}(a)
  =i\circ\hat{g}\circ g^{\ast}(a)=i\circ g(a)=g(a).
\end{align*}
Hence we conclude that $f\circ h=g$. Therefore $\boldsymbol{R}$ is onto\=/projective.
\end{proof}

In order to show that our abstract framework is well behaved, we are committed to proving that every equivalence between two substitution\=/invariant \mdr's is induced by a pair of endomoprhisms on the $\boldsymbol{\Sigma}_{\mathcal{L}}$\=/module $\Mult_\mathcal{L}$ (note that this claim can be seen as a variant of the Isomorphism Theorem in the setting of \mdr's). In other words, we want to show that $\Mult_{\mathcal{L}}$ has the REP (that is, it is onto\=/projective) in the category of $\boldsymbol{\Sigma}_{\mathcal{L}}$\=/modules. Instead of proving this directly, we will take a brief detour and prove some more general results. First, we make a note of the following definition.

\begin{definition}
An $\mathbf{A}$\=/module $\boldsymbol{R}$ is \emph{cyclic} if there is $a\in R$ such that $R=\{\sigma\ast a : \sigma\in A\}$.
\end{definition}

Observe that every dually integral po\=/semiring $\mathbf{A} = \langle A,\leq,+,\cdot,0,1\rangle $ can be seen as a degenerate instance of $\mathbf{A}$\=/module if we drop $\cdot$ and $1$ from the signature and set $\ast=\cdot$. Keeping this in mind, we obtain the following:

\begin{lemma}\label{Lem:Envelope}
Any dually integral po\=/semiring $\mathbf{A}$, viewed as an $\mathbf{A}$\=/module, is cyclic and onto\=/projective.
\end{lemma}

\begin{proof}
Clearly $\mathbf{A}$ is cyclic, since $A=\{\sigma\cdot1 : \sigma\in A\}$. Let $f\colon\boldsymbol{R}\to\boldsymbol{S}$ and $g\colon\mathbf{A}\to\boldsymbol{S}$ be two morphisms, where $f$ is onto. Fix any $a\in R$ such that $f(a)=g(1)$ and define $h\colon A\to R$ via $h(\sigma)  =\sigma\ast^{\boldsymbol{R}}a$. Clearly, $h\colon\mathbf{A}\to\boldsymbol{R}$ is a morphism. Moreover, given $\sigma\in A$, we have that
\[
 g(\sigma) = g(\sigma\cdot1) = \sigma\cdot g(1) = \sigma\cdot f(a)
 = f(\sigma\ast a) = f \circ h(\sigma).
\]
Hence we conclude that $\mathbf{A}$ is onto\=/projective.
\end{proof}

Cyclic modules can be described in an arrow-theoretic way as follows:

\begin{lemma}\label{Lem:Cyclic}
An $\mathbf{A}$\=/module $\boldsymbol{R}$ is cyclic iff there is an onto morphism $f\colon\mathbf{A}\to\boldsymbol{R}$.
\end{lemma}

\begin{proof}
If there is an onto morphism $f\colon\mathbf{A}\to\boldsymbol{R}$ and $x\in R$, then for some $\sigma\in A$ we have that $x=f(\sigma) =\sigma\ast f(1)  $. To prove the converse, it is enough to check that if $R=\{\sigma\ast v : \sigma\in A\}$, then the map $f\colon\mathbf{A} \to\boldsymbol{R}$ defined by $f(\sigma)  =\sigma\ast v$ is a morphism.
\end{proof}

We are now ready to prove the following characterisation of cyclic and onto\=/projective objects in $\Mod{\mathbf{A}}$.

\begin{theorem}\label{Thm:CyclicProjective}
Let $\boldsymbol{R}$ be an $\mathbf{A}$\=/module. The following conditions are equivalent:
\begin{enumerate}
\item $\boldsymbol{R}$ is cyclic and onto\=/projective.

\item There is a retraction $f\colon\mathbf{A}\to\boldsymbol{R}$.

\item There are $\mu\in A$ and $v\in R$ such that $\mu\ast v=v$ and $A\ast\{v\} = R$ and for every $\sigma,\pi\in A$: if $\sigma\ast v\leq\pi\ast v$, then $\sigma\cdot\mu\leq\pi\cdot\mu$.
\end{enumerate}
\end{theorem}

\begin{proof}
To prove 1.\ implies 2.\ observe that from Lemma~\ref{Lem:Cyclic}, we know that there is a morphism $f\colon\mathbf{A}\to\boldsymbol{R}$ which is surjective. Applying the projectivity of $\boldsymbol{R}$ to the diagram given by $f$ and the identity map $id_{\boldsymbol{R}}$, we conclude that $f$ is a retraction.

Next we prove that 2.\ implies 1. From Lemma~\ref{Lem:Cyclic} we know that $\boldsymbol{R}$ is cyclic. Moreover, $\boldsymbol{R}$ is a retract of an onto\=/projective object by Lemma~\ref{Lem:Envelope}. Thus we conclude that $\boldsymbol{R}$ is onto\=/projective too.

To prove 2.\ implies 3.\ observe that by the assumption, there is an injective morphism
$g\colon\boldsymbol{R}\to\mathbf{A}$ such that $1_{\boldsymbol{R}} = f\circ g$. Then we define $v=f(1)$ and $\mu=g(v)$. Since $f$ is onto, we have that $A\ast\{v\}=R$. Moreover:
\[
 \mu\ast v = \mu\ast f(1) = f(\mu\ast1) = f(\mu) = f(g(v)) = v.
\]
Considering $\sigma,\pi\in A$ such that $\sigma\ast v\leq\pi\ast v$, we have that $\sigma\cdot\mu=\sigma\cdot g(v)=g(\sigma\ast v)\leq g(\pi\ast v)=\pi\cdot g(v)=\pi\cdot\mu.$

\newpage

Finally we prove 3. implies 2. Since $A\ast\{v\}=R$, we know that the map $f\colon \mathbf{A}\to\boldsymbol{R}$ defined as $f(\sigma) = \sigma\ast v$ is an onto morphism. Then let $g\colon\boldsymbol{R} \to\mathbf{A}$ be defined via $g(\sigma\ast v) =\sigma\cdot\mu$. Using the assumption, it is not difficult to see that $g$ is well defined and order-preserving. Also, it can be routinely established that $g$ preserves the action and is a monoid homomorphism. Thus, $g$ is a morphism. In order to prove that $f\circ g=1_{\boldsymbol{R}}$ we consider a generic element $\sigma\ast v\in R$ and show that:
\[
 f\circ g(\sigma\ast v) = \sigma\ast(f\circ g(v)) = \sigma\ast f(\mu)
 = \sigma\ast(\mu\ast v) = \sigma\ast v.
\]
\end{proof}

\begin{theorem}\label{blascoferrer}
The $\boldsymbol{\Sigma}_{\mathcal{L}}$\=/module $\Mult_{\mathcal{L}}$ is cyclic and onto\=/projective. In particular, this implies that it has the REP.
\end{theorem}

\begin{proof}
Let $x$ be a designated $\mathcal{L}$-variable, and let $v=[x]$. Moreover, let $\sigma$ be the $\mathcal{L}$-substitution defined by $\sigma(y) = x$ for all $\mathcal{L}$-variables $y$, and fix $\mu=[\sigma]  $. It is not difficult to see that $v$ and $\mu$ satisfy the conditions of Item (3) in Theorem~\ref{Thm:CyclicProjective}.
\end{proof}

\section{Multiset deductive relations}\label{energade}

Recall that what prompted us to extend Blok and J\'{o}nsson's theory was the motivating example of \emph{multiset deductive relations} (\mdr's), defined in Definition~\ref{d:MDR}.
It turns out that our general theory has interesting offshoots once we focus on this special case --- and the whole of the present section will be devoted to buttressing this claim. For a start, we list some prototypical instances of \mdr's.

\begin{example}\label{commu}
Recall that an algebra $\mathbf{A}=\langle A,\wedge,\vee ,\cdot,\to,1\rangle$ of type $\mathcal{L}_{0}=\langle 2,2,2,2,0\rangle$ is a \emph{commutative and integral residuated lattice} (see e.g.~\cite{MPT}) if $\langle A,\wedge,\vee\rangle$ is a lattice, $\langle A,\cdot,1\rangle$ is a commutative monoid, $1$ is the top element w.r.t.\ the induced order $\leq$ of $\langle A,\wedge,\vee\rangle$, and the following residuation law holds for every $a,b,c\in A$:
\[
 a\cdot b\leq c \quad\iff\quad a\leq b\to c.
\]
Given a class $\mathcal{K}$ of commutative and integral residuated lattices, let the relation $\vdash_{\mathcal{K}}$ be defined as follows for all $\Gamma = [\varphi_{1},\dots,\varphi_{n}]  ,\Delta=[\psi _{1},\dots,\psi_{m}]  \in \Fm_{\mathcal{L}_{0}}^{\flat}$:\footnote{Here and in the sequel, given a multiset $\Gamma = [\varphi_{1},\dots,\varphi _{n}]  $ of $\mathcal{L}_{0}$\=/formulas, the notation $\varphi_{1}\cdot\ldots\cdot\varphi_{n}$ will ambiguously refer to any of the $\mathcal{L}_{0}$\=/formulas
\[
  (\cdots(\varphi_{f(1)}\cdot\varphi_{f(2)})  \cdot \ldots \cdot\varphi_{f(n)}),
\]
where $f$ is a permutation of $\{1,\dots,n\}  $. By way of
convention, if $\Gamma$ is the empty multiset, we formally set $\varphi_{1}\cdot\ldots\cdot\varphi_{n}=1$.}
\begin{equation}\label{Eq:Consequence}
 \Gamma\vdash_{\mathcal{K}}\Delta \quad\iff\quad
 \mathcal{K}\models \varphi_{1}\cdot\ldots\cdot\varphi_{n}\leq\psi_{1}\cdot\ldots\cdot\psi_{m}.
\end{equation}
It can be checked that $\vdash_{\mathcal{K}}$ is a substitution\=/invariant \mdr\ in the sense of Definition~\ref{d:MDR}.
\end{example}

Example~\ref{commu} identifies, for every substructural logic whose equivalent algebraic semantics is a quasi-variety of commutative and integral residuated lattices, a ``multiset-theoretic'' companion of such that best suits the resource interpretation at which we hinted in our introduction. One particular such logic will play some role in what follows. The multiset companion $\vdash_{\mathcal{MV}}$ of infinite-valued \L{}ukasiewicz logic $\vdash_{\text{\L}}$ is obtained when the class $\mathcal{K}$ is the variety $\mathcal{MV}$ of \emph{MV-algebras} (see \cite{CDM}), formulated in the language $\mathcal{L}_{0}$.\footnote{On the other hand, the logic obtained when $\mathcal{K}$ is the variety generated by the $3$-element MV-chain was first considered by Arnon Avron (see e.g.~\cite{Av94}).}

Also, observe that Example~\ref{commu} encompasses the so-called \emph{internal consequence relations} of algebraisable substructural sequent calculi with exchange and weakening (see~\cite{Avronl, Avron}). In fact, let $S$ be such a calculus and $\mathcal{Q}$ its equivalent algebraic semantics. Upon defining, for finite multisets of $\mathcal{L}_{0}$\=/formulas $\Gamma$ and $\Delta=[\psi_{1},\dots,\psi_{m}]$,
\[
 \Gamma\vdash_{\mathcal{S}}\Delta \quad\iff\quad
 \vdash_{S}\Gamma \to \psi_{1}\cdot\ldots\cdot\psi_{m},
\]
then it follows from well\=/known results about substructural logics that ${\vdash_{\mathcal{S}}} = {\vdash_{\mathcal{Q}}}$.

If the above examples look a bit contrived, this is due, in part, to the fact that the multiple-conclusion format is unfamiliar to many. As a consequence, it would seem expedient to extract from these examples appropriate single\=/conclusion relations that can be more easily compared, say, to the usual, \emph{external}, consequence relations of substructural sequent calculi. It turns out that single\=/conclusion relations can be recovered as \emph{fragments} of \mdr's (see~\cite{CP}):

\begin{definition}\label{d:consequenc}
Let $\mathcal{L}$ be a propositional language. A \emph{single\=/conclusion} \mdr\ on $\mathcal{L}$ is a relation ${\vdash^{u}}\subseteq \Fm_{\mathcal{L}}^{\flat}\times \Fm_{\mathcal{L}}$ such that, for some \mdr\ $\vdash$,
\[
 \Gamma\vdash^{u}\alpha \quad\iff\quad \Gamma\vdash[\alpha].
\]
\end{definition}

It should be observed, though, that we are not claiming that single\=/conclusion \mdr's be themselves instances of \mdr's, for they need not be closed w.r.t.\ all the conditions that define them (see~\cite{CP} for further discussion). Clearly, for each single\=/conclusion \mdr\ ${\vdash^{u}}$ there exists the least \mdr\ ${\vdash}$ that has ${\vdash^{u}}$ as fragment: namely, the intersection of all such \mdr's. In Subsection~\ref{ss:Hilbert} we present an example of two \mdr's with the same single\=/conclusion fragment.

Multiset deductive relations can be taken to subsume \tcr's, as illustrated by the next example.

\begin{example}
Every finitary substitution\=/invariant \tcr\ can be encoded into a finitary substitution\=/invariant \mdr. Indeed, consider such \tcr\ $\Vdash$ on language $\mathcal{L}$. Then we define a substitution\=/invariant \mdr\ $\vdash$ on $\mathcal{L}$ by setting, for all $\Gamma = [\varphi_{1},\dots,\varphi_{n}]  ,\Delta=[\psi_{1},\dots,\psi_{m}]  $ in $\Fm_{\mathcal{L}}^{\flat}$:
\[
 \Gamma\vdash\Delta
 \quad\iff\quad
 \vert \Gamma\vert \Vdash \psi_{k}, \text{ for all }k\leq m.
\]
It is a purely computational matter to check that $\vdash$ is indeed a substitution\=/invariant \mdr. Moreover, it is clear that $\vdash$ encodes $\Vdash$ in the sense that, whenever $\varphi_{j}\neq\varphi_{k}$ for all $j,k\leq n$,
\[
 \varphi_{1},\dots,\varphi_{n}\Vdash\psi
 \quad\iff\quad
 [\varphi_{1},\dots,\varphi_{n}] \vdash [\psi].
\]
\end{example}

\subsection{Hypermatrices and the first completeness theorem}

In this subsection we describe a matrix-based semantics for arbitrary substitution\=/invariant \mdr's. To this end, we work in a fixed (but otherwise arbitrary) language $\mathcal{L}$.

Logical matrices are part and parcel of every algebraic logician's toolbox (see e.g.\ \cite[Ch.~4]{Font}). As a consequence, when we are dealing with \mdr's over a language $\mathcal{L}$, it seems desirable to be in a position to help ourselves to concepts that inherit at least some of the effectiveness and power of matrix semantics in AAL. Whatever notion of matrix we are bound to adopt, it appears natural that its attendant notion of ``Lindenbaum--Tarski matrix'' be in keeping with Definition~\ref{teoria}: we expect such matrices to have the form $\langle \mathbf{Fm}_{\mathcal{L}},F\rangle $ where $F$ is a certain \emph{set of finite multisets} of $\mathcal{L}$\=/formulas. Therefore, it is all too plausible to focus on ``matrices'' constituted by an algebra and a certain family of finite submultisets of its universe. The next definitions spell out in detail this basic insight.

\begin{definition}\label{HyperMatrices}
An $\mathcal{L}$-\emph{hypermatrix} is a pair $\langle\mathbf{A},F\rangle$, where $\alg{A}$ is an $\mathcal{L}$-algebra and $F$ a $\leqslant$-downset in $A^{\flat}$.
\end{definition}

\begin{definition}\label{HyperMatrices1}
For a class $\mathsf{H}$ of $\mathcal{L}$-hypermatrices we define a relation $\models_{\Mat{H}}$ on $\Fm_{\mathcal{L}}^{\flat}$ as $\Gamma \models_{\Mat{H}}\Delta$, if for every $\mat{A}=\langle\alg{A},F\rangle\in\Mat{H}$, each context $\mathfrak{C}\in A^{\flat}$ and each homomorphism $f\colon\mathbf{Fm}_{\mathcal{L}}\to\alg{A}$:
\[
 \mathfrak{C}\uplus f(\Gamma)\in F
 \quad\then\quad
 \mathfrak{C}\uplus f(\Delta)\in F.
\]
\end{definition}

\begin{theorem}\label{beddamatri}
Let $\Mat{H}$ be a class of $\mathcal{L}$-hypermatrices. Then $\models_{\Mat{H}}$ is a substitution\=/invariant \mdr\ on $\mathcal{L}$.
\end{theorem}

\begin{proof}
We show the proof for $\Mat{H}=\{\langle\alg{A},F\rangle\}$; the general statement then follows from the obvious facts that ${\models_{\Mat{H}}} = \bigcap \{\models_{\mat{A}} : \mat{A}\in\Mat{H}\}  $ and that the class of substitution\=/invariant \mdr's is closed under intersections. The validity of transitivity and substitution\=/invariance of $\models_{\mat{A}}$ is very easy to see.

For Compatibility, assume that $\Gamma\models_{\mat{A}}\Delta$ and consider a  context ${\mathfrak{C}}\in {A^{\flat}}$ and a homomorphism $e\colon\mathbf{Fm}_{\mathcal{L}}\to\alg{A}$. If ${\mathfrak{C}}\uplus e(\Gamma)\uplus e(\Pi)={\mathfrak{C}}\uplus e(\Gamma\uplus\Pi)\in F$, then by our hypothesis ${\mathfrak{C}}\uplus e(\Delta)\uplus e(\Pi)\in F$, and thus $\Gamma\uplus \Pi\models_{\mat{A}}\Delta{\uplus\Pi}$.

For Generalised Reflexivity, assume that $\Gamma\leqslant\Delta$ and consider a context ${\mathfrak{C}}\in {A^{\flat}}$ and a homomorphism $e\colon\mathbf{Fm}_{\mathcal{L}}\to\alg{A}$. Clearly, $e(\Gamma)  \leqslant e(\Delta)  $, and so by the compatibility of $\leqslant$ we obtain ${\mathfrak{C}\uplus e}(\Gamma) \leqslant{\mathfrak{C}\uplus e}(\Delta)  $. Thus, if ${\mathfrak{C}\uplus e}(\Delta)  \in F$, then ${\mathfrak{C} \uplus e}(\Gamma)  $, because $F$ is a $\leqslant$-downset.
\end{proof}

Note that the fact that $F$ a $\leqslant$-downset in $A^{\flat}$, and the reference to arbitrary contexts $\mathfrak{C}$, play a crucial role in the previous proof. Lifting both restrictions at once leads to the following definition of a relation $\models_{\mat{A}}'$ on $\Fm_{\mathcal{L}}^{\flat}$ for an arbitrary pair $\mat{A}=\langle\alg{A},F\rangle$, where $\alg{A}$ is an $\mathcal{L}$-algebra and $F\subseteq A^{\flat}$:
\[
 \Gamma\models_{\mat{A}}'\Delta
 \text{ iff for\ each homomorphism} f\colon\mathbf{Fm}_{\mathcal{L}}\to\alg{A}, \text{ if } f(\Gamma)\in F\text{ then }f(\Delta)\in F.
\]

\begin{lemma}\label{l:WeirdConsequence}
Consider a pair $\mat{A} = \langle\alg{A},F\rangle$, where $\alg{A}$ is an $\mathcal{L}$-algebra and $F\subseteq A^{\flat}$. Then ${\models_{\mat{A}}}\subseteq{\models_{\mat{A}}'}$. Assume further that for some \mdr\ $\vdash$ we have ${\vdash}\subseteq{\models_{\mat{A}}'}$. Then $\mat{A}$ is an $\mathcal{L}$-hypermatrix and\/ ${\vdash}\subseteq{\models_{\mat{A}}}$.
\end{lemma}

\begin{proof}
The first inclusion is trivial. Assume that ${\vdash}\subseteq{\models_{\mat{A}}'}$ and we show that $F$ is $\leqslant$-downset. Note that for any multiset of mutually different atoms $[p_{1},\dots,p_{n}]$ and $m\leq n$ we have $[p_{1},\dots,p_{n}]\models'_{\mat{A}}[p_{1},\dots,p_{m}]$ and for any multisets $\mathfrak{X}\leqslant\mathfrak{Y}$ there is homomorphism $e\colon\mathbf{Fm}_{\mathcal{L}}\to\alg{A}$  such that $\mathfrak{X}=[e(p_{1}),\dots,e(p_{m})]$ and $\mathfrak{Y}=[e(p_{1}),\dots,e(p_{n})]$.

To complete the proof we need to show that ${\vdash}\subseteq{\models_{\mat{A}}}$. Assume that $\mathfrak{C}=[x_{1},\dots,x_{n}]\in A^{\flat}$, $\Gamma\vdash\Delta$  and $f\colon\mathbf{Fm}_{\mathcal{L}} \to\alg{A}$ are such that $\mathfrak{C}\uplus f(\Gamma)\in F$ and we need to prove that $\mathfrak{C}\uplus f(\Delta)\in F$. Let $\Pi$ be a multiset of mutually different atoms $[p_{1},\dots,p_{n}]$ not occurring in $\Gamma\uplus\Delta$; we know that $\Gamma\uplus\Pi\vdash\Delta\uplus\Pi$. Next consider the homomorphism $e'$ defined as $e'(p_{i})=x_{i}$ and $e'(p)=e(p)$ for other atoms and note that $e'(\Gamma\uplus\Pi)=\mathfrak{X} \uplus e(\Gamma)\in F$ and so $\mathfrak{X}\uplus e(\Gamma)=e' (\Delta\uplus\Pi)\in F$.
\end{proof}

Next, we provide an example showing that $\models_\mat{A}$ and $\models'_\mat{A}$ are in general different relations.

\begin{example}
Consider a two-element set $A = \{0, 1 \}$. Then let $F \subseteq A^{\flat}$ be defined as follows:
\[
 F = \{\emptyset, [0], [1], [0, 1] \}.
\]
Now, equip $A$ with the structure of an algebra $\mathbf{A} = \langle A, \bold{0}, \bold{1}\rangle$, whose only operations are constant symbols $\bold{0}$ and $\bold{1}$ for $0$ and $1$, respectively. Clearly, $F$ is a $\leqslant$-downset and it is easy to see that for $\mat{A} = \langle\alg{A},F\rangle$ we have
\[
 [\bold{0}] \models_{\mat{A}}' [\bold{1}]
 \quad\text{and}\quad
 [\bold{0}, \bold{1}] \not\models_{\mat{A}}'[\bold{1}, \bold{1}].
\]
Hence the consequence $\models_{\mat{A}}'$ does not satisfy Compatibility and, therefore, it is not an \mdr\ and cannot be equal to $\models_\mat{A}$.
\end{example}

\begin{corollary}
Let $\Mat{H}$ be a class of $\mathcal{L}$-hypermatrices. We define a relation $\models_{\Mat{H}}'$ on $\Fm_{\mathcal{L}}^{\flat}$ as ${\models_{\Mat{H}}'} = \bigcap\{\models'_{\mat{A}} : \mat{A}\in\Mat{H}\}$. Then $\models_{\Mat{H}}'$ is an \mdr\ iff\/ ${\models_{\Mat{H}}'}={\models_{\Mat{H}}}$.
\end{corollary}

Now we can define notions of model and filter. Note that the previous lemma renders it immaterial whether we use ${\models_{\mat{A}}'}$ or ${\models_{\mat{A}}}$ in such definitions.

\begin{definition}\label{Eq:FirstModels}
Let\/ $\vdash$ be a substitution\=/invariant \mdr\ on $\mathcal{L}$. An $\mathcal{L}$\=/hypermatrix $\mat{A}=\langle\alg{A},F\rangle$ is a \emph{model} of\/ $\vdash$ and $F$ is an $\vdash$\=/\emph{filter} on $\alg{A}$ if\/ ${\vdash}\subseteq{\models_{\mat{A}}}$. By\/ $\mathbf{Mod}({\vdash })  $ we denote the set of all models of\/ ${\vdash}$ and by $\mathcal{F}i_{{\vdash}}(\alg{A})  $ the set of all\/ ${\vdash}$-filters on~$\alg{A}$.
\end{definition}

It is straightforward to show that $\mathcal{F}i_{\vdash}(\alg{A})$ is a closure system. Given a closure system $\mathcal{C}$ on a set $X$, let us denote by $\mathcal{C}^{p}$ the set of its principal members:
\[
\mathcal{C}^{p} = \Bigl\{\bigcap\{C\in\mathcal C : x\in C\} : x\in X \Bigr\}.
\]

\begin{proposition}\label{p:FiltersAreDS}
Let $\vdash$ be a substitution\=/invariant \mdr\ on $\mathcal{L}$.
\begin{enumerate}
\item For every $\mathcal{L}$-algebra $\alg{A}$, the collection $(\mathcal{F}i_{{\vdash}}(\alg{A})  )^{p}$ is a \ds\ on $\mathbf{A}^{\flat}$ and moreover $\mathcal{F}i_{{\vdash}}(\alg{A}) = ((\mathcal{F}i_{{\vdash}}(\alg{A})  )^{p})^{\mathit{BJ}}$.

\item $(\mathcal{F}i_{{\vdash}}(\mathbf{Fm}_{\mathcal{L}})  )^{p}=Th^p(\vdash)$.
\end{enumerate}
\end{proposition}

\begin{proof}
Consider an arbitrary $\mathcal{L}$-algebra $\mathbf{B}$. Recall that $\mathcal{F}i_{\vdash }(\mathbf{B})  $ is a closure system. Then let $\textup{Fg}_{\vdash}^{\alg{B}}\colon\wp(B^{\flat})\to\wp(B^{\flat})$ be its corresponding closure operator. It is easy to see that for every $X\subseteq B^{\flat}$,
\[
 \textup{Fg}_{\vdash}^{\alg{B}}(X) = \bigcup_{n\in\omega}F_{n}
\]
where $F_{0}=X$ and
\begin{align*}
 F_{n+1} = F_{n}\cup\{\mathfrak{X}\in B^{\flat} :\
  &\text{there are } \Gamma,\Delta\in \Fm_{\mathcal{L}}^{\flat}\text{ s.t.\ }\Gamma\vdash\Delta \text{ and a homomorphism}\\
  & f\colon\mathbf{Fm}_{\mathcal{L}}\to\mathbf{B}\text{ s.t.\ }f(\Gamma)\in F_{n} \text{ and } f(\Delta) = \mathfrak{X}\}.
\end{align*}
From the above remarks it follows that for every $\Gamma\in \Fm_{\mathcal{L}}^{\flat}$,
\[
 \textup{Fg}_{\vdash}^{\mathbf{Fm}_{\mathcal{L}}}(\Gamma) = \{\Delta\in \Fm_{\mathcal{L}}^{\flat} : \Gamma\vdash\Delta\}.
\]
In particular, this means that $(\mathcal{F}i_{{\vdash}}(\mathbf{Fm}_{\mathcal{L}}))^{p}=Th^{p}(\vdash)$, which proves the second statement.

For the other statement, consider an $\mathcal{L}$-algebra $\alg{A}$. We begin by proving that $(\mathcal{F}i_{{\vdash}}(\alg{A}))^{p}$ is a \ds\ on $\langle A^{\flat},\leq,\uplus,\emptyset\rangle$. To this end, we claim that \begin{equation}\label{Eq:Claim} \text{if }  \mathfrak{X}\in\textup{Fg}_{\vdash}^{\alg{A}}(\mathfrak{Y}), \text{ then }\mathfrak{X}\uplus\mathfrak{C}\in\textup{Fg}_{\vdash}^{\alg{A}}(\mathfrak{Y}\uplus\mathfrak{C})
\end{equation}
for every $\mathfrak{X},\mathfrak{Y},\mathfrak{C}\in A^{\flat}$.

To prove this claim, fix $\mathfrak{X},\mathfrak{Y},\mathfrak{C} = [c_{1},\dots,c_{k}]\in A^{\flat}$ and consider decompositions
\[
 \textup{Fg}_{\vdash}^{\alg{A}}(\mathfrak{Y}) = \bigcup_{n\in\omega}F_{n}
 \quad\text{and}\quad
 \textup{Fg}_{\vdash}^{\alg{A}}(\mathfrak{Y}\uplus\mathfrak{C}) = \bigcup_{n\in\omega}G_{n}
\]
defined at the beginning of this proof. We show, by induction on $n\in\omega$, that
\begin{equation}\label{Eq:InductiveTrick}
 \mathfrak{X}\in F_{n} \quad\then\quad \mathfrak{X}\uplus\mathfrak{C}\in G_{n}.
\end{equation}
The case where $n=0$ is direct. Then we consider the case $n=s+1$. Suppose that $\mathfrak{X}\in F_{s+1}$. Then there are $\Gamma,\Delta$ and a homomorphism $f\colon\mathbf{Fm}_{\mathcal{L}}\to\alg{A}$ such that
\[
 \Gamma\vdash\Delta,\ f(\Gamma)\in F_{s}, \text{ and } f(\Delta)= \mathfrak{X}.
\]
Consider the multiset $\Pi = [x_1,\dots,x_k]$ consisting of fresh pairwise different variables. By compatibility of $\vdash$, we have that
\begin{equation}\label{Eq:Filter}
 \Gamma\uplus\Pi\vdash\Delta \uplus \Pi.
\end{equation}
By inductive hypothesis we know that
\begin{equation}\label{Eq:Induction}
 f(\Gamma)\uplus \mathfrak{C} \in G_{s}.
\end{equation}
Let $f'\colon\mathbf{Fm}_{\mathcal{L}}\to\alg{A}$ be any homomorphism which coincides with $f$ on the variables appearing in the formulas $\Gamma$ and $\Delta$, and such that $f'(x_{i})=c_{i}$. By~\eqref{Eq:Induction} we have that
\[
 f'(\Gamma)\uplus f'(\Pi)\in G_{s}.
\]
Together with~\eqref{Eq:Filter}, this implies that
\[
 \mathfrak{X}\uplus\mathfrak{C} = f(\Delta)\uplus \mathfrak{C} = f'(\Delta)\uplus f'(\Pi) \in G_{s+1}.
\]
This concludes the proof of~\eqref{Eq:InductiveTrick} and, therefore, establishes~\eqref{Eq:Claim}. Now we turn back to the main argument. First observe that
\[
 (\mathcal{F}i_{{\vdash}}(\alg{A}))^{p} = \{\textup{Fg}_{\vdash}^{\alg{A}}(\mathfrak{X}) : \mathfrak{X}\in A^{\flat}\}.
\]
Clearly $(\mathcal{F}i_{{\vdash}}(\alg{A}))^{p}$ is a family of $\leqslant$-downsets. Define the map $\delta\colon A^{\flat}\to(\mathcal{F}i_{{\vdash}}(\alg{A}))^{p}$ setting
\[
 \delta(\mathfrak{X}) = \bigcap\{C\in(\mathcal{F}i_{{\vdash}}(\alg{A}))^{p} : \mathfrak{X}\in C\}
 =\textup{Fg}_{\vdash}^{\alg{A}}(\mathfrak{X})
\]
for every $\mathfrak{X}\in A^{\flat}$. Clearly $\delta(A^{\flat})=(\mathcal{F}i_{{\vdash}}(\alg{A}))^{p}$. Finally, from~\eqref{Eq:Claim} we obtain that
\[
 \delta(\mathfrak{X})\subseteq\delta(\mathfrak{Y}) \quad\then\quad
 \delta(\mathfrak{X}\uplus\mathfrak{C})\subseteq\delta(\mathfrak{Y}\uplus\mathfrak{C}).
\]
Hence we conclude that $(\mathcal{F}i_{{\vdash}}(\alg{A}))^{p}$ is a \ds\ on $\mathbf{A}^{\flat}$ as desired.

Then we turn to prove that $\mathcal{F}i_{{\vdash}}(\alg{A})=((\mathcal{F}i_{{\vdash}}(\alg{A}))^{p})^{\mathit{BJ}}$. The inclusion from left to right is clear. To prove the other inclusion, consider a family $\{F_{i} : {i\in I}\}\subseteq(\mathcal{F} i_{{\vdash}}(\alg{A}))^{p}$. Then suppose that $\Gamma\vdash\Delta$, and consider a homomorphism $f\colon\mathbf{Fm}_{\mathcal{L}}\to\alg{A}$ such that $f(\Gamma)\in\bigcup_{i\in I}F_{i}$. Clearly there is $j\in I$ such that $f(\Gamma)\in F_{j}$. Since $F_{j}\in\mathcal{F}i_{\vdash}(\alg{A})$, we obtain that $f(\Delta)\in F_{j}\subseteq\bigcup_{i\in I}F_{i}$. Hence we conclude that $\bigcup_{i\in I}F_{i}\in\mathcal{F}i_{\vdash}(\alg{A})$.
\end{proof}

The notions introduced so far are enough to obtain a first completeness theorem for any substitution\=/invariant \mdr.

\begin{theorem}[1st Completeness Theorem]
Let ${\vdash}$ be a substitution\=/invariant \mdr\ on $\mathcal{L}$. Then
\[
{\vdash} = {\models_{\mathbf{Mod}({\vdash})}} = {\models'_{\mathbf{Mod}({\vdash})}}.
\]
\end{theorem}

\begin{proof}
From left to right, our claim is obvious. For the reverse direction, assume that $\Gamma\not \vdash \Delta$ and define $T=\mathrm{Th}_{{\vdash}}(\mathrm{\Gamma})$. By Proposition~\ref{p:FiltersAreDS}, $\langle \mathbf{Fm}_{\mathcal{L}},T\rangle  \in\mathbf{Mod}({\vdash})$ and then the identity mapping is the homomorphism we need to show that $\Gamma\not \models_{\mathbf{Mod}({\vdash})}\Delta$.
\end{proof}

\subsection{A bridge to Gentzen systems and the second completeness theorem}

We will now establish a connection with the algebraic theory of Gentzen systems, i.e.\ substitution\=/invariant \acr's on sequents, a well\=/trodden research stream in AAL (\cite{Ra06,JRa13,ReV93,ReV95,Py99,T91}), that will serve as a touchstone for our approach based on hypermatrices. Let $\vdash$ be a substitution\=/invariant \mdr\ on $\mathcal{L}$. We will associate with it a consequence relation $\vdash^{g}$ between sequents. To this end, consider the set $Seq_{\mathcal{L}}$ of $\mathcal{L}$\=/sequents of the form
\[
 \emptyset \rhd \langle\varphi_{1},\dots,\varphi_{n}\rangle,
\]
where $\langle\varphi_{1},\dots,\varphi_{n}\rangle$ is a finite sequence of formulas and the relation ${\vdash^{g}}\subseteq {\mathcal{\wp}(Seq_{\mathcal{L}})\times Seq_{\mathcal{L}}}$ defined as follows:
\begin{align*}
 X\vdash^{g}\emptyset \rhd\langle\varphi_{1},\dots,\varphi_{n}\rangle
 \quad\iff\quad
 &\text{ there is }\emptyset \rhd\langle\gamma_{1},\dots,\gamma_{m}\rangle\in X\\
 &\text{ s.t.\ } [\gamma_{1},\dots,\gamma_{m}]\vdash [\varphi_{1},\dots,\varphi_{n}].
\end{align*}

Clearly, $\vdash^{g}$ is a substitution\=/invariant \acr\ on $Seq_{\mathcal{L}}$: for every substitution $\sigma$, if $X\vdash^{g}\emptyset \rhd \langle\varphi_{1},\dots,\varphi_{n}\rangle$, then
\[
 \{\emptyset\rhd\langle\sigma(\gamma_{1}),\dots,\sigma(\gamma_{m})\rangle : \emptyset\rhd\langle\gamma_{1},\dots,\gamma_{m}\rangle\in X\}\vdash^{g}\emptyset\rhd\langle\sigma(\varphi_{1}),\dots,\sigma(\varphi_{n})\rangle.
\]

As we remarked above, substitution\=/invariant \acr's on sequents are the object of study of numerous papers that have appeared under the heading of \emph{algebraisation of Gentzen systems}. Within this theory, a model of a Gentzen system $\Vdash$ on $Seq_{\mathcal{L}}$ is a pair $\langle\alg{A},F\rangle$ where $\alg{A}$ is an $\mathcal{L}$-algebra and $F$ is a set of finite sequences of elements of $A$ such that for every set $X\cup\{\emptyset \rhd\langle\varphi_{1},\dots,\varphi_{n}\rangle\}\subseteq Seq_{\mathcal{L}}$, if $X\Vdash\emptyset  \rhd\langle\varphi_{1},\dots,\varphi_{n}\rangle$, then for every homomorphism $f\colon\mathbf{Fm}_{\mathcal{L}}\to\alg{A}$,
\begin{multline}\label{Eq:SecondModels}
\text{if }\langle f(\gamma_{1}),\dots,f(\gamma_{m})\rangle\in F  
\text{ for every }\emptyset\rhd\langle\gamma_{1},\dots,\gamma_{m}\rangle\in X,\\
 \text{ then }\langle f(\varphi_{1}),\dots,f(\varphi_{n})\rangle\in F.
\end{multline}
We denote by $\mathbf{Mod}(\Vdash)$ the class of all models of $\Vdash$. A quick comparison between Definition~\ref{Eq:FirstModels} (see also the comments before the definition) and~\eqref{Eq:SecondModels} suggests that the models of $\vdash$ and $\vdash^{g}$ must be interdefinable. To make this idea precise, we define two maps as follows:
\begin{equation}\label{Eq:Transformations}
 (\cdot)^{s}\colon\mathbf{Mod}(\vdash) \longleftrightarrow \mathbf{Mod} (\vdash^{g}) \cocolon(\cdot)^{m}
\end{equation}
where $(\cdot)^{s}$ stands for  \emph{sequents} and $(\cdot)^{m}$ stands for  \emph{multisets}. For $\langle\alg{A},F\rangle\in{\mathbf{Mod} (\vdash^{g})}$ and $\langle\mathbf{B},G\rangle\in\mathbf{Mod}(\vdash)$, we set:
\begin{align*}
  \langle\alg{A},F\rangle^{m}
  &= \langle\alg{A},\{[ a_{1},\dots,a_{n}] : n\geq 0, \langle a_{1},\dots,a_{n}\rangle\in F\}\rangle\\
  \langle\mathbf{B},G\rangle^{s}
  &= \langle\mathbf{B},\{f\in B^{\{1,\dots,n\}  } :  n\geq 0, [f(1), \dots, f(n)]=\mathfrak{X} \text{ for some } \mathfrak{X} \in G\}\rangle.
\end{align*}

The proof of the following result is simple (note that we need to use Lemma~\ref{l:WeirdConsequence}):

\begin{lemma}\label{Lem:FirstCorrespondence}
The transformations $(\cdot)^{s} \colon\mathbf{Mod}(\vdash)\longleftrightarrow\mathbf{Mod}(\vdash^{g}) \cocolon (\cdot)^{m}$ are well defined and mutually inverse bijections.
\end{lemma}

As a consequence, we can apply the algebraic constructions developed for Gentzen systems in the above-mentioned literature, to the study of substitution\=/invariant \mdr's. We devote the remaining part of this subsection to give a flavour of the resulting theory.

Let $\langle\alg{A},F\rangle$ be a pair consisting of an $\mathcal{L}$-algebra $\alg{A}$ and a set $F$ of finite sequences of elements of $A$. A congruence $\theta$ of $\alg{A}$ is  \emph{compatible} with $F$ if for every $a_{1},b_{1},\dots,a_{n},b_{n}\in A$,
\[
\text{if }\langle a_{1},\dots,a_{n}\rangle\in F\text{ and } \langle a_{1},b_{1}\rangle,\dots,\langle a_{n},b_{n}\rangle\in\theta,\text{ then }\langle b_{1},\dots,b_{n}\rangle\in F.
\]
When $\theta$ is compatible with $F$, we set
\[
F/\theta=\{\langle a_{1}/\theta,\dots,a_{n}/\theta\rangle : \langle a_{1},\dots,a_{n}\rangle\in F\}.
\]
It turns out that there exists the largest congruence of $\alg{A}$ compatible with $F$. This congruence is called the \emph{Leibniz congruence} of $F$ over $\alg{A}$, and is denoted by $\Omega^{\alg{A}}F$. The  \emph{reduced models} of $\Vdash$ are the following class:
\[
 \mathbf{Mod}^{\ast}(\Vdash) =
 \mathbb{I} \bigl\{\langle\alg{A}/\Omega^{\alg{A}}F,F/\Omega^{\alg{A}}F\rangle :
 \langle\alg{A},F\rangle\in\mathbf{Mod}(\Vdash) \bigr\}.
\]
A general result~\cite[Proposition~5.111]{Font} shows that the Gentzen system $\Vdash$ is  \emph{complete} with respect to the semantics $\mathbf{Mod}^{\ast}(\Vdash)$: the pairs in $\mathbf{Mod}^{\ast}(\Vdash)$ are models of $\Vdash$ and, moreover, if $X\nVdash\emptyset \rhd\langle \varphi_{1},\dots,\varphi_{n}\rangle$, then there is $\langle\alg{A},F\rangle\in\mathbf{Mod}^{\ast}(\Vdash)$ and a homomorphism $f\colon \mathbf{Fm}_{\mathcal{L}}\to\alg{A}$ such that
\begin{multline}
\langle f(\gamma_{1}),\dots,f(\gamma_{m})\rangle\in F  
\text{ for every } \emptyset\rhd\langle\gamma_{1},\dots,\gamma_{m}\rangle\in X,\\
\text{ and }\langle f(\varphi_{1}),\dots,f(\varphi_{n})\rangle\notin F.
\end{multline}

All the above constructions can be transferred to the study of substitution\=/invariant \mdr's as follows. Let $\langle\alg{A},F\rangle$ be an $\mathcal{L}$-hypermatrix. A congruence $\theta$ of $\alg{A}$ is  \emph{compatible} with $F$ if for every $a_{1},b_{1},\dots,a_{n},b_{n}\in A$,
\begin{multline}
\text{if }[  a_{1},\dots,a_{n}]  \in F  
\text{ and } [a_{1}/\theta,\dots,a_{n}/\theta] = [b_{1}/\theta,\dots,b_{n}/\theta],\\
\text{ then }[  b_{1},\dots,b_{n}] \in F.
\end{multline}
When $\theta$ is compatible with $F$, we set:
\[
 F/\theta = \{[a_{1}/\theta,\dots,a_{n}/\theta]  : [a_{1},\dots,a_{n}]  \in F\}.
\]

\begin{lemma}\label{Lem:ExistenceOfLeibniz}
There exists the largest congruence of $\alg{A}$ compatible with $F$.
\end{lemma}

\begin{proof}
We will denote by $G$ the set of finite sequences of elements of $A$ such that $\langle\alg{A},G\rangle=\langle\alg{A},F\rangle^{s}$. We know that $\Omega^{\alg{A}}G$ is the largest congruence of $\alg{A}$ compatible with $G$. In order to conclude the proof, it will be enough to show that $\Omega^{\alg{A}}G$ is also the largest congruence of $\alg{A}$ compatible with $F$. However, it is easily proved that a congruence of $\alg{A}$ is compatible with $F$ if and only if it is compatible with $G$, whence our claim follows.
\end{proof}

Given the above result, we denote the largest congruence of $\alg{A}$ compatible with $F$ by $\Omega^{\alg{A}}F$, and call it the the  \emph{Leibniz congruence} of $F$ over $\alg{A}$. We define the \emph{reduced models} of an \mdr\ $\vdash$ as follows:
\[
 \mathbf{Mod}^{\ast}(\vdash) =
 \mathbb{I} \bigl\{ \langle\alg{A}/\Omega^{\alg{A}}F, F/\Omega^{\alg{A}}F\rangle :
 \langle\alg{A},F\rangle\in\mathbf{Mod}(\vdash) \bigr\}.
\]

\begin{corollary}
The transformations $(\cdot)^{s}\colon\mathbf{Mod}^{\ast}(\vdash) \longleftrightarrow \mathbf{Mod}^{\ast}(\vdash^{g})\cocolon(\cdot)^{m}$ are well defined and mutually inverse bijections.
\end{corollary}

\begin{proof}
Pick a model $\langle\alg{A},F\rangle\in\mathbf{Mod}^{\ast}(\vdash)$. Let $G$ be the set of finite sequences of elements of $A$ such that $\langle\alg{A},F\rangle^{s}=\langle\alg{A},G\rangle$. From Lemma~\ref{Lem:FirstCorrespondence} we know that $\langle\alg{A},G\rangle \in\mathbf{Mod}(\vdash^{g})$. Moreover, since $\langle\alg{A},F\rangle \in\mathbf{Mod}^{\ast}(\vdash)$, the congruence $\Omega^{\alg{A}}F$ is the identity relation $\textup{Id}_{\alg{A}}$. Now, in the proof of Lemma~\ref{Lem:ExistenceOfLeibniz} we showed that $\Omega^{\alg{A}}F = \Omega^{\alg{A}}G$. Thus we conclude that $\Omega^{\alg{A}}G=\textup{Id}_{\alg{A}}$. Since $\langle\alg{A},G\rangle\cong\langle\alg{A}/\textup{Id}_{\alg{A}},G/\textup{Id}_{\alg{A}}\rangle$ and $\mathbf{Mod}^{\ast}(\vdash^{g})$ is closed under isomorphisms, we conclude that $\langle\alg{A},F\rangle^{s}=\langle\alg{A},G\rangle \in\mathbf{Mod}^{\ast}(\vdash^{g})$. A similar argument shows that if $\langle\mathbf{B},G\rangle\in\mathbf{Mod}^{\ast}(\vdash^{g})$, then $\langle\mathbf{B},G\rangle^{m}\in\mathbf{Mod}^{\ast}(\vdash)$. This means that the maps $(\cdot)^{s}\colon\mathbf{Mod}^{\ast}(\vdash) \longleftrightarrow \mathbf{Mod}^{\ast}(\vdash^{g})\cocolon(\cdot)^{m}$ are well defined. The fact that they are inverse bijections follows from Lemma~\ref{Lem:FirstCorrespondence}.
\end{proof}

\begin{theorem}[2nd completeness theorem]\label{t:CompletenessReduced}
Let ${\vdash}$ be a substitution\=/invariant \mdr\ on $\mathcal{L}$. Then
\[
 {\vdash} = {\models_{\mathbf{Mod}^\ast({\vdash})}} = {\models'_{\mathbf{Mod}^\ast({\vdash})}}.
\]
\end{theorem}

\begin{proof}
From the general theory of the algebraisation of Gentzen systems we know that $\mathbf{Mod}^{\ast}(\vdash ^{g})$ is a class of models of $\vdash^{g}$. By Lemma~\ref{Lem:FirstCorrespondence} this implies that the image of the class $\mathbf{Mod}^{\ast}(\vdash^{g})$ under the transformation $(\cdot)^{m}$ is a class of models of $\vdash$. But by the previous corollary we know that this image coincides with $\mathbf{Mod}^{\ast}(\vdash)$. Thus $\mathbf{Mod}^{\ast}(\vdash)$ is a class of models of $\vdash$. Then suppose that $\mathfrak{X}  \nvdash \mathfrak{Y}$. From the very definition of $\vdash^{g}$ it follows that for any sequence $\langle\varphi_{1},\dots,\varphi_{n}\rangle$ and $\langle\psi_{1},\dots,\psi_{m}\rangle$ such that $\mathfrak{X} = [\varphi_{1},\dots,\varphi_{n}]$ and $\mathfrak{Y}=[\psi_{1},\dots,\psi_{m}]$ we have:
\[
 \emptyset\rhd\langle\varphi_{1},\dots,\varphi_{n}\rangle \nVdash\emptyset
 \rhd\langle\psi_{1},\dots,\psi_{m}\rangle.
\]
From the general completeness result of $\vdash^{g}$ with respect to $\mathbf{Mod}^{\ast}(\vdash^{g})$ it follows that there is $\langle \alg{A},F\rangle\in\mathbf{Mod}^{\ast}(\vdash^{g})$ and a homomorphism $f\colon\mathbf{Fm}_{\mathcal{L}}\to\alg{A}$ such that
\[
 \langle f(\varphi_{1}),\dots,f(\varphi_{n})\rangle\in F
 \quad\text{and}\quad
 \langle f(\psi_{1}),\dots,f(\psi_{m})\rangle\notin F.
\]
Let $G\subseteq A^{\flat}$ be such that $\langle\alg{A},F\rangle^{m}=\langle\alg{A},G\rangle$. Due to the previous corollary we know that $\langle\alg{A},G\rangle\in\mathbf{Mod}^{\ast}(\vdash)$. Moreover, it is straightforward to check that
\[
 [f(\varphi_{1}),\dots,f(\varphi_{n})]  \in G
 \quad\text{and}\quad
 [f(\psi_{1}),\dots,f(\psi_{m})] \notin G.
\]
\end{proof}

\subsection{Monoid matrices and t-norm semantics}

It is possible to give hypermatrices a (nearly) equivalent formulation in such a way as to shed further light on the direction in which they generalise ordinary logical matrices. The rough idea is replacing the unstructured set of designated values in a logical matrix by a richer structure. If $\langle \alg{A},D\rangle $ is an ordinary logical matrix (i.e., an algebra with a subset), $D$ can be identified with a function in $\{0,1\}  ^{A}$; in other words, being designated is an all-or-nothing matter. The set $\{0,1\}  $ can also be viewed as the universe of the $2$-element join semilattice $\mathbf{2}$. If we replace $\mathbf{2}$ by any dually integral Abelian pomonoid $\mathbf{D}$, however, we can at the same time express an ordering of ``degrees of designation'', and evaluate the degree of designation of whole submultisets of $A$, with the monoidal operation in $\mathbf{D}$ ensuring that evaluations behave well with respect to multiset union. This leads to the following:

\begin{definition}
Let $\mathcal{L}$ be a language. An $\mathcal{L}$-\emph{monoid matrix} is a quadruple $\mat{M}=\langle \alg{A},\mathbf{D},G,f\rangle $, where:
\begin{enumerate}
\item $\alg{A}$ is an $\mathcal{L}$-algebra.

\item $\mathbf{D} = \langle D,\leq,+,0\rangle$ is a dually integral Abelian pomonoid.

\item $G$ is a $\leq$-downset in $\mathbf{D}$.

\item $f\colon \mathbf{A}^\flat \to \langle D,\leq,+,0\rangle $ is a pomonoid homomorphism.
\end{enumerate}
\end{definition}

The next lemma ensures that monoid matrices are closed w.r.t.\ a sort of quotient construction. For $E\subseteq D$, $(E]  $ will denote the $\leq$-downset generated in $\mathbf{D}$ by $E$.

\begin{lemma}
Let $\mat{M} = \langle \alg{A},\mathbf{D},G,f\rangle $ be an $\mathcal{L}$-monoid matrix, $\mathbf{D}'$  a dually integral Abelian pomonoid, and  $g\colon\mathbf{D}\to\mathbf{D}'$ a pomonoid homomorphism. Then
\[
g_{\mathbf{D}'}(\mat{M}) = \langle \alg{A},\mathbf{D}',(g(G)],g\circ f\rangle
\]
is an $\mathcal{L}$-monoid matrix.
\end{lemma}

\begin{proof}
$\mathbf{D}'$ is a dually integral Abelian pomonoid by assumption, and likewise $(g(G)]$ is a $\leq$-downset in $\mathbf{D}'$. The map $g\circ f$ is a composition of monoid homomorphisms, and if $\mathfrak{X} \leqslant\mathfrak{Y}$, then $f(\mathfrak{X}) \leq^{\mathbf{D}}f(\mathfrak{Y})$ and so $g(f(\mathfrak{X})) \leq^{\mathbf{D}'}g(f(\mathfrak{Y}))$.
\end{proof}

The relationships between the previously introduced notions are made clear in
the next theorem. While every $\mathcal{L}$-hypermatrix arises out
of an $\mathcal{L}$-monoid matrix, an $\mathcal{L}$-monoid matrix need not be
more than the ``homomorphic image'' of an $\mathcal{L}$-monoid matrix that
arises out of an $\mathcal{L}$-hypermatrix.

\begin{theorem}\mbox{}
\begin{enumerate}
\item If $\mat{M}=\langle \alg{A},\mathbf{D},G,f\rangle$ is an $\mathcal{L}$-monoid matrix, then
\[
 \mat{H}^{\mat{M}} = \langle \alg{A},\{\mathfrak{X}\in A^{\flat} : f(\mathfrak{X})\in G\}\rangle
\]
is an $\mathcal{L}$-hypermatrix.

\item If $\mat{H} = \langle\alg{A},F\rangle$ is an $\mathcal{L}$-hypermatrix, then
\[
\mat{M}^{\mat{H}} = \langle \alg{A}, \alg{A}^\flat,F,id\rangle
\]
is an $\mathcal{L}$-monoid matrix.

\item $\mat{H}^{\mat{M}^{\mat{H}}}=\mat{H}$.

\item $f_{\mathbf{D}}(\mat{M}^{\mat{H}^{\mat{M}}})
=\mat{M}$.
\end{enumerate}
\end{theorem}

\begin{proof}
To prove the first claim we have to show that $\{\mathfrak{X}\in A^{\flat} : f(\mathfrak{X})  \in G\}  $ is a $\leqslant$-downset, i.e.\ if $f(\mathfrak{X}) \in G$ and $\mathfrak{Y} \leqslant\mathfrak{X}$, then $f(\mathfrak{Y})  \in G$. However, if $\mathfrak{Y} \leqslant\mathfrak{X}$, then $f(\mathfrak{Y})  \leq^{\mathbf{D}}f(\mathfrak{X})$, whence our conclusion follows as $G$ is a $\leq^{\mathbf{D}}$-downset.

The second claim is trivial and the third one as straightforward:
\[
\mat{H}^{\mat{M}^{\mat{H}}} = \langle \alg{A},\{\mathfrak{X}\in A^{\flat} : id(\mathfrak{X}) \in F\} \rangle =\langle \alg{A},F\rangle = \mat{H}.
\]

The final claim: $\mat{M}^{\mat{H}^{\mat{M}}}=\langle \alg{A},\alg{A}^\flat,\{\mathfrak{X} : f(\mathfrak{X})  \in G\},id\rangle $, whence
\[
 f_{\mathbf{D}}\Big(\mat{M}^{\mat{H}^{\mat{M}}}\Big)
 = \langle\alg{A},\mathbf{D},(f(\{\mathfrak{X} : f(\mathfrak{X})  \in G\})],f\rangle.
\]
However, $f(\{\mathfrak{X} : f(\mathfrak{X}) \in G\}) \subseteq G$, so $(f(\{\mathfrak{X} : f(\mathfrak{X})  \in G\})] = G$, and thus $f_{\mathbf{D}} \Big(\mat{M}^{\mat{H}^{\mat{M}}}\Big) = \mat{M}$.
\end{proof}

We now focus on a special class of monoid matrices, namely, those matrices whose underlying pomonoid is just the closed unit real interval $[0,1]$, endowed with some t-norm $\ast$ (i.e., a monotone, associative, and commutative operation with unit 1) and with the usual ordering of real numbers. In essence, these monoid matrices can be seen as an algebra together with a fuzzy set of designated values. Although very special in nature, these matrices can be used to yield a semantics for the multiset companions of some fuzzy logics, obtained in the same guise as $\vdash_{\mathcal{MV}}$ (see the remarks immediately following Example~\ref{commu}).

\begin{definition}
An $\mathcal{L}$-\emph{fuzzy matrix} is an $\mathcal{L}$-monoid matrix $\mat{M} = \langle \alg{A},\mathbf{D},G,f\rangle$, s.t.\ $\mathbf{D} = \langle [  0,1]  ,\sqsupseteq,\ast,1\rangle $, where $\ast$ is some t-norm and $\sqsubseteq$ is the usual ordering of $[  0,1]  $.
\end{definition}

Two fuzzy matrices will be called \emph{similar} if their algebra reducts are similar and the t-norm $\ast$ is the same in both cases.

\begin{definition}
If $\mat{M} = \langle \alg{A},\mathbf{D},G,f\rangle$ is an $\mathcal{L}$-fuzzy matrix with t-norm $\ast$ and $\Gamma,\Delta\in \Fm_{\mathcal{L}}^{\flat}$, we set $\Gamma\models_{\mat{M}}^{\ast}\Delta$ just in case $\Gamma\models_{\mat{H}^{\mat{M}}}\Delta$.

If $\Mat{M}$ is a class of similar $\mathcal{L}$-fuzzy matrices, we write $\Gamma\models_{\Mat{M}}^{\ast}\Delta$ as a shortcut for: $\Gamma \models_{\mat{M}}^{\ast}\Delta$ for every $\mat{M}\in\Mat{M}$.
\end{definition}

Theorem~\ref{beddamatri} implies the following:

\begin{lemma}\label{lem:FuzzyMatrices}
If\/ $\Mat{M}$ is a class of similar $\mathcal{L}$-fuzzy matrices, $\models_{\Mat{M}}^{\ast}$ is a substitution\=/invariant \mdr\ on $\mathcal{L}$.
\end{lemma}

Observe that, if we fix $\alg{A}$ and $f$ while letting $\Mat{M}$ be the class of similar $\mathcal{L}$-fuzzy matrices $\{\langle \alg{A},\mathbf{D}, [a,1], f\rangle : a\in[0,1]\}$, we have that
\begin{align*}
 \Gamma\models_{\Mat{M}}^{\ast}\Delta
  &\iff \Gamma \models_{\mat{H}^{\mat{M}}}\Delta \text{ for all } \mat{M}\in\Mat{M}\\
  &\iff \forall a\forall\mathfrak{C}\forall e (f (\mathfrak{C}\uplus e(\Gamma)) \in[a,1] \then f (\mathfrak{C}\uplus e(\Delta)\in[a,1]))\\
  &\iff \forall\mathfrak{C}\forall e(f(\mathfrak{C} \uplus e(\Gamma)) \sqsubseteq f(\mathfrak{C} \uplus e(\Delta)))\\
  &\:\then \forall e(f(e(\Gamma))\sqsubseteq f(e(\Delta))),
\end{align*}
where the third equivalence uses the downward closure of
$\mat{H}^{\mat{M}}$. Furthermore, since
\[
 f(e[\gamma_{1},\dots,\gamma_{n}]) = f(e(\gamma_{1}))\ast\cdots\ast f(e(\gamma_{n})),
\]
the behaviour of $f$ is entirely determined by its behaviour on one-element multisets, whence we lose no generality in taking $f$ to be a function from $A$ to $[0,1]$. In other words, a fuzzy matrix can be viewed --- in this special case --- as an algebra together with a fuzzy set of designated values. Moreover, if $\alg{A}$ itself is some algebra with universe $[0,1]$, the function $f$ becomes a real function.

With this material at hand, we are ready to prove the completeness of $\vdash_{\mathcal{MV}}$ with respect to the class $\Mat{M}$ of all $\mathcal{L}_{0}$-fuzzy matrices of the form $\langle [\mathbf{0,1}]_{\mathrm{MV}},\mathbf{D},[a,1], f\rangle$, where:
\begin{enumerate}
\item[$\bullet$] $[\mathbf{0,1}]_{\mathrm{MV}}$ is the standard $\mathrm{MV}$ algebra over $[0,1]$, formulated in the language $\mathcal{L}_{0}$ of commutative residuated lattices.

\item[$\bullet$] $\mathbf{D}=\langle [0,1],\sqsupseteq,\otimes,1\rangle $, where $\otimes$ is the \L{}ukasiewicz t-norm.\footnote{When discussing \L{}ukasiewicz logic and its multiset companion $\vdash_{\mathcal{MV}}$, we write the multiplicative conjunction (residuated lattice product) $\varphi\cdot\psi$ using the more customary notation $\varphi\otimes\psi$.}

\item[$\bullet$] $a\in[0,1]$.

\item[$\bullet$] $f\colon[0,1] \to [0,1]$ is strictly monotone and preserves $\otimes$ (this class is non\=/empty: it contains e.g.\ the identity function $id$ and the square function $(\cdot)^{2}$).
\end{enumerate}

\begin{theorem}\label{matricks}
For any $\Gamma, \Delta\in\Fm_{\mathcal{L}_{0}}^{\flat}$, the following are equivalent:
\begin{enumerate}
\item $\Gamma\vdash_{\mathcal{MV}}\Delta$.

\item $\Gamma\models_{\Mat{M}}^{\ast}\Delta$.

\item $\Gamma\models_{\Mat{I}}^{\ast}\Delta$, where $\Mat{I}=\{\mat{M}\in\Mat{M}:  f=id\}$.
\end{enumerate}
\end{theorem}

\begin{proof}
We prove that 1. implies 2., then next implication is trivial and the final one
follows from the observation after  Lemma~\ref{lem:FuzzyMatrices}

Suppose  that $\Gamma=[  \varphi_{1},\dots,\varphi_{n}]$, $\Delta=[\psi_{1},\dots,\psi_{m}]$, $\Gamma\vdash_{\mathcal{MV}}\Delta$. By Chang's completeness theorem, this means that for every homomorphism $e\colon \mathbf{Fm}_{\mathcal{L}_{0}}\to[  \mathbf{0,1}]_{\mathrm{MV}}$,%
\[
 e(\varphi_{1})  \otimes\cdots\otimes e(\varphi_{n}) \sqsubseteq e(\psi_{1}) \otimes\cdots\otimes e(\psi_{m}).
\]

Let $\mat{M}=\langle [  \mathbf{0,1}]_{\mathrm{MV}},\mathbf{D},[a,1]  ,f\rangle \in\Mat{M}$, $\mathfrak{C\in}[\mathbf{0,1}]  ^{\flat}$, and let $e'\colon\mathbf{Fm}_{\mathcal{L}_{0}}\to[  \mathbf{0,1}]_{\mathrm{MV}}$ be a homomorphism. Suppose further that
\[
a    \sqsubseteq f(\mathfrak{C\uplus}e'(\Gamma))
  = f(\mathfrak{C})  \otimes f(e'(\Gamma)  )
  = f(\mathfrak{C})  \otimes f(e'(\varphi_{1})  )  \otimes\cdots\otimes f(e'(\varphi_{n})).
\]
However, since $f$ is monotone and $\otimes$-preserving,
\[
 f(e'(\varphi_{1})) \otimes\cdots\otimes f(e'(\varphi_{n})) \sqsubseteq f(e'(\psi_{1}))  \otimes\cdots\otimes f(e'(\psi_{m}))
\]
and, by monotonicity of t-norms, $a\sqsubseteq f(\mathfrak{C})\otimes f(e'(\varphi_{1})) \otimes \cdots\otimes f(e'(\varphi_{n})) \sqsubseteq f(\mathfrak{C})  \otimes f(e'(\psi_{1})  ) $ $ \otimes \cdots\otimes f(e'(\psi_{m})  ) =f(\mathfrak{C\uplus}e'(\Delta))$, which suffices for our conclusion.
\end{proof}

\subsection{Hilbert systems}\label{ss:Hilbert}

We mentioned at the outset that previous attempts at investigating multiset consequence are few and far between. Virtually all authors who undertook this enterprise, however, tried to set up axiomatic calculi of sorts (\cite{Avron, Troelstra, Primer}). We now proceed to present our own take on the issue.

\begin{definition}\label{consecution}
A \emph{consecution}\footnote{The term ``consecution'' is taken from~\cite{AB} (the term ``sequent'' is sometimes used instead).} in a propositional language $\mathcal{L}$ is a pair $\langle  {\Gamma,\Delta }\rangle  $, where $\Gamma$ and $\Delta$ are finite multisets of formulas. A consecution is \emph{single\=/conclusion} if $\Delta=[\varphi]$ for some formula $\varphi$.
\end{definition}

Instead of `$\langle  {\Gamma,\Delta}\rangle$', we write `$\Gamma  \rhd\Delta$'. With a slight abuse, we also identify the consecution $\emptyset \rhd\lbrack\varphi]$ with the formula
$\varphi$.

\begin{definition}[Axiomatic system]\label{d:AxSys}
Let $\mathcal{L}$ be a propositional language. A \emph{(single-conclusion) axiomatic system} in the language $\mathcal{L}$ is a set $\mathsf{AS}$ of (single\=/conclusion) consecutions closed under arbitrary substitutions.\footnote{I.e., if $\Gamma \rhd\Delta \in \mathsf{AS}$, then $\sigma(\Gamma) \rhd\sigma(\Delta) \in \mathsf{AS}$.} The elements of $\mathsf{AS}$ of the form $\Gamma \rhd\Delta$ are called \emph{axioms} if\/ $\Gamma=\emptyset$ and \emph{deduction rules} otherwise.
\end{definition}

Of course, each axiomatic system can also be seen as a collection of \emph{schemata}, i.e.\ a collection of consecutions and all their substitution instances. Observe that our single\=/conclusion axiomatic systems are essentially Avron's \emph{multiset Hilbert systems} (see~\cite{Avron}); however, the upcoming notion of tree-proof is different, as our single\=/conclusion \mdr's (unlike Avron's ``simple consequence relations'') enjoy the Monotonicity condition.

\begin{definition}[Tree-proof]\label{d:pruvu}
Let $\mathcal{L}$ be a propositional language and let $\mathsf{AS}$ be a single\=/conclusion axiomatic system in $\mathcal{L}$. A \emph{tree-proof} of a formula $\varphi$ from a multiset of formulas $\Gamma$ in $\mathsf{AS}$ is a finite tree $t$ labelled by formulas such that:
\begin{enumerate}
\item[$\bullet$] The root of $t$ is labelled by $\varphi$.

\item[$\bullet$] If a leaf of $t$ is labelled by $\psi$, then either
\begin{enumerate}
\item[$-$] $\psi$ is an axiom or
\item[$-$] $\psi$ is an element of $\Gamma$ and it labels at most $\Gamma(\psi)$ leaves in $t$.
\end{enumerate}

\item[$\bullet$] If a node of $t$ is labelled by $\psi$ and $\Delta\neq\emptyset$ is the multiset of labels of its predecessor nodes, then $\Delta \rhd \lbrack\psi]\in\mathsf{AS}$.
\end{enumerate}
\end{definition}

We write $\Gamma\vdash_{\!\mathsf{AS}}^{t}\varphi$ whenever there is a tree-proof of $\varphi$ from $\Gamma$ in $\mathsf{AS}$. Our next goal is to define a notion of derivation for arbitrary axiomatic systems.

\begin{definition}\label{d:proof}
Let ${\mathcal{L}}$ be a propositional language and let $\mathsf{AS}$ be an axiomatic system in ${\mathcal{L}}$. A \emph{derivation} of a finite multiset of formulas $\Delta$ from a finite multiset of formulas $\Gamma$ in $\mathsf{AS}$ is a finite sequence $\langle \Gamma_{1},\dots,\Gamma_{n}\rangle$ of finite multisets of formulas such that:
\begin{enumerate}
\item[$\bullet$] $\Gamma_{1}=\Gamma$;

\item[$\bullet$] For every $\Gamma_{j}$, $1<j\leq n$, there is a rule $\Psi \rhd \Psi'\in\mathsf{AS}$, such that $\Psi\leqslant\Gamma_{j-1}$ and $\Gamma_{j}=(\Gamma_{j-i}\setminus\Psi)\uplus\Psi'$;

\item[$\bullet$] $\Delta\leqslant\Gamma_{n}$.
\end{enumerate}

We say that $\Delta$ is \emph{derivable} from $\Gamma$ in $\mathsf{AS}$, and
write $\Gamma\vdash_{\!\mathsf{AS}}\Delta$, if there is a derivation of\/
$\Delta$ from $\Gamma$ in $\mathsf{AS}$.
\end{definition}

Observe that, if $\Psi=\emptyset$, the second clause above says that{ in a derivation} we are allowed to beef up with finitely many axioms any multiset that has already been derived. The next lemma supports the adequacy of the given definition.

\begin{lemma}\label{prelu}
Let ${\mathcal{L}}$ be a propositional language and let $\mathsf{AS}$ be an axiomatic system in ${\mathcal{L}}$. Then $\vdash_{\!\mathsf{AS}}$ is the least substitution\=/invariant \mdr\ containing $\mathsf{AS}$.
\end{lemma}

\begin{proof}
Generalised Reflexivity being trivial, we prove the remaining conditions.
\begin{description}
\item[Compatibility] Given a derivation $P=\langle \Gamma_{1},\dots,\Gamma_{n}\rangle $ of $\Delta$ from $\Gamma$ in $\mathsf{AS}$, it is easy to observe that the sequence $P'=\langle \Gamma_{1}\uplus\Pi,\dots,\Gamma_{n} \uplus\Pi\rangle $ is a derivation of $\Delta\uplus\Pi$ from $\Gamma\uplus\Pi$ (thanks to the fact that in our notion of proof we can apply rules in an arbitrary context).

\item[Transitivity] Suppose we have a derivation $\langle \Gamma_{1} ,\dots,\Gamma_{n}\rangle $ of $\Delta$ from $\Gamma$ in $\mathsf{AS}$. Then $\Delta\leqslant\Gamma_{n}$ and so from $\Delta\vdash_{\!\mathsf{AS}}\Pi$ we get by monotony $\Gamma_{n}\vdash_{\!\mathsf{AS}}\Pi\uplus(\Gamma _{n}\setminus\Delta)$; let $\langle \Delta_{1},\dots,\Delta_{n} \rangle $ be the corresponding derivation in $\mathsf{AS}$. Note that $\Delta_{1}=\Gamma_{n}$ and $\Pi\leqslant\Delta_{n}$. Then the sequence
\[
\langle \Gamma_{1},\dots,\Gamma_{n},\Delta_{2},\dots,\Delta_{n}\rangle
\]
is clearly a derivation of $\Pi$ from $\Gamma$ in $\mathsf{AS}$.

\item[Substitution\=/invariance] Given a derivation $P=\langle \Gamma _{1},\dots,\Gamma_{n}\rangle $ of $\Delta$ from $\Gamma$ in $\mathsf{AS}$, the sequence $P'=\langle \sigma[\Gamma_{1}] ,\dots,\sigma[\Gamma_{n}]  \rangle $ is a derivation of $\sigma[\Delta]$ from $\sigma[\Gamma]$ in $\mathsf{AS}$.
\end{description}

Now for the proof that $\vdash_{\!\mathsf{AS}}$ is the least substitution\=/invariant \mdr\ containing $\mathsf{AS}$. Obviously $\mathsf{AS} \subseteq{\vdash_{\mathsf{AS}}}$. What remains to prove is that for each substitution\=/invariant \mdr\ $\vdash$, if $\mathsf{AS}\subseteq\vdash$, then ${\vdash_{\mathsf{AS}}}\subseteq\,\vdash$. Assume that $\Gamma\vdash _{\mathsf{AS}}\Delta$, i.e.\ there is a derivation $P$ of $\Delta$ from $\Gamma$ in $\mathsf{AS}$. By induction on the length of $P$, we can show that for each multiset of formulas $\Pi$ in $P$ we have $\Gamma\vdash\Pi$, and hence in particular $\Gamma\vdash\Delta$. The base case is settled with an appeal to Reflexivity. As to the induction step: let $\Pi$ and $\Pi'$ be labels of successive elements of $P$ and $\Gamma\vdash\Pi$. We know that there is a rule $\Psi \rhd\Psi'$ such that $\Pi' = (\Pi\setminus\Psi)\uplus\Psi'$. Thus $\Psi\vdash\Psi'$, and so by Compatibility $\Psi\uplus(\Pi\setminus\Psi)  \vdash \Psi'\uplus(\Pi\setminus\Psi)  $, i.e.\ $\Pi\vdash \Pi'$. An application of Transitivity completes the proof.
\end{proof}

\begin{definition}\label{d:presentation}
Let ${\mathcal{L}}$ be a propositional language, $\mathsf{AS}$ an axiomatic system in~${\mathcal{L}}$, and let\/ $\vdash$ be a substitution\=/invariant \mdr\ on~${\mathcal{L}}$. We say that $\mathsf{AS}$ is an \emph{axiomatic system} for (or a \emph{presentation} of) $\vdash$ if\/ ${\vdash} = {\vdash_{\mathsf{AS}}}$.
\end{definition}

Clearly, due to the previous lemma, each \mdr\ can be seen as its own presentation, and so we obtain:

\begin{corollary}[\L{}os--Suszko]\label{szko}
Every given substitution\=/invariant \mdr\ $\vdash$ coincides with the derivability relation $\vdash_{\mathsf{AS}}$ of some axiomatic system $\mathsf{AS}$.
\end{corollary}

The next lemma spells out the relation between derivations and tree-proofs.

\begin{lemma}\label{l:trees}
Let $\mathsf{AS}$ be a single\=/conclusion axiomatic system on~${\mathcal{L}}$, $\Gamma\vdash_{\!{\mathsf{AS}}}\Delta$, and $\varphi \in\vert \Delta\vert $. Then there are multisets of~${\mathcal{L}} $\=/formulas $\Gamma^{\varphi}$ and $\Gamma^{r}$ such that $\Gamma^{\varphi }\uplus\Gamma^{r}=\Gamma$, $\Gamma^{\varphi}\vdash_{\!{\mathsf{AS}}} ^{t}\varphi$ and $\Gamma^{r}\vdash_{\!{\mathsf{AS}}}\Delta\setminus \lbrack\varphi]$.
\end{lemma}

\begin{proof}
Let $P=\langle \Gamma_{1},\dots,\Gamma_{n}\rangle $ be the assumed derivation of $\Delta$ from $\Gamma$ in ${\mathsf{AS}}$. For each $i\leq n$, let $\Gamma_{i}=[\psi_{1}^{i},\dots,\psi_{k_{i}}^{i}]$ and note that without loss of generality we can assume  that the rule used in the $i$-th step of $P$ is $[\psi_{1}^{i}, \dots, \psi_{p_i}^{i} ] \rhd\psi_{c_{i}}^{i+1}$ for some $p_i \leq k_{i}$ and $c_{i}\leq k_{i+1}$. Also note that $k_{i+1}-1=k_{i+1}- p_{i}$ and there is a bijection $f$ between $\Gamma_{i+1}\setminus\lbrack\psi_{c_{i}}^{i+1}]$ and $\Gamma _{i}\setminus[\psi_{1}^{i}, \dots, \psi_{p_i}^{i} ]$ such that $\psi_{j} ^{i+1}=\psi_{f(j)  }^{i}$ whenever $c_{i}\neq j\leq k_{i+1}$. We construct the labelled graph $G$ with nodes $N=\{\langle i,j\rangle\mid i\leq n\text{ and }j\leq k_{i}\}$, where $\psi_{j}^{i}$ is the label of $\langle i,j\rangle$, and edges only between the following nodes:
\begin{enumerate}
\item[$\bullet$] \emph{rule edges}: $\langle i,k\rangle$ and $\langle{i+1},c_{i}\rangle$
for each $k\leq p_{i};$

\item[$\bullet$] \emph{non-rule edges}: $\langle i,f(j)\rangle$ and $\langle{i+1},j\rangle$ for $j\neq c_{i}$.
\end{enumerate}

It is easy to see that $G$ is a forest (a disjoint union of trees). Let $t$ be the subtree of $G$ with root $\psi^{n}_{1}$, and let $\Gamma^{\varphi}$ be the multiset of all labels of leaves in $t$ which are not axioms. Then clearly $\Gamma^{\varphi}\leqslant\Gamma$ and $t$ is \emph{almost} a tree-proof of $\varphi$ from $\Gamma^{\varphi}$; all we have to do is to collapse nodes connected by non-rule edges.

Finally, let $\Gamma_{i}^{t}$ denote the multiset resulting from $\Gamma_{i}$ by removing  formulas labeling the nodes of $t$ (as many times as it labels some node) and observe that $\Gamma_{1}^{t},\dots,\Gamma_{n}^{t}$ is \emph{almost} a proof of $\Gamma_{n}^{t}=\Gamma_{n}\setminus\lbrack\psi_{1}^{n}]$ from $\Gamma_{1}^{t}$: we only need to remove each $\Gamma_{i}^{t}$ which equals its predecessor. Defining $\Gamma^{r}=\Gamma_{1}^{t}$ and observing that $\Gamma^{r} =\Gamma\setminus\Gamma^{\varphi}$ and $\Delta\setminus\lbrack\psi_{1} ^{n}]\leqslant\Gamma_{n}^{t}$ completes the proof.
\end{proof}

\begin{lemma}
For any $\mathsf{AS}$ single\=/conclusion axiomatic system, $\Gamma \vdash_{\!{\mathsf{AS}}}^{t}\varphi$ iff\/ $\Gamma\vdash_{\!{\mathsf{AS}}}[\varphi]$.
\end{lemma}

\begin{proof}
One direction follows directly from the previous lemma. To prove the converse one, assume that there is a tree-proof $t$ of $\varphi$ from $\Gamma$ in $\mathsf{AS}$. Let $n$ be a node of $t$, $\psi_{n}$ its label, $P_{n}$ the set of its predecessors, $\Delta_{n}$ the multiset of labels of nodes in $P_{n}$, and $\Gamma_{n}$ the multiset of labels of elements of $\Gamma$ which are not axioms and occur in leaves of the subtree of $t$ with root $n$. If we show that $\Gamma_{n}\vdash_{\!\mathsf{AS}}[\psi_{n}]$, the proof is done: indeed for the root $r$ of $t$ we obtain $\Gamma_{r}\vdash_{\!\mathsf{AS}}[\varphi]$ and given that $\Gamma_{r} \leqslant\Gamma$ and $\vdash$ is an \mdr, we obtain the claim by Monotonicity. Let us prove the claim: if $n$ is a leaf, the proof is trivial. Otherwise, there is a rule $\Delta_{n} \rhd\lbrack \psi_{n}]$ and for each $m\in P_{n}$ we have $\Gamma_{m}\vdash\lbrack\psi _{m}]$. Thus $\biguplus\limits_{m\in P_{n}}\Gamma_{m}\vdash\Delta_{n}$ and so $\biguplus\limits_{m\in P_{n}}\Gamma_{m}\vdash\lbrack\psi_{n}]$. The proof is completed by observing that $\biguplus\limits_{m\in P_{n}}\Gamma_{m} =\Gamma_{n}$.
\end{proof}

Recall that single\=/conclusion \mdr's were introduced in Definition~\ref{d:consequenc} as \emph{fragments} of \mdr's. We now observe that the tree-provability relations of single\=/conclusion axiomatic systems, in a sense, ``generate'' the corresponding derivability relations.

\begin{corollary}
Let $\mathsf{AS}$ be a single\=/conclusion axiomatic system. Then\/ $\vdash ^{t}_{\!{\mathsf{AS}}}$ is a single\=/conclusion \mdr\ and\/ $\vdash_{\!{\mathsf{AS} }}$ is the least\/ \mdr\ $\vdash$ such that
\[
 \Gamma\vdash^{t}_{\!{\mathsf{AS}}}\alpha
 \quad\iff\quad
 \Gamma\vdash[\alpha].
\]
\end{corollary}

\begin{proof}
The former claim is immediate from the previous lemma. To prove the latter, assume that $\vdash$ is an \mdr\ such that
\[
 \Gamma\vdash_{\!{\mathsf{AS}}}^{t}\alpha
 \quad\iff\quad
 \Gamma\vdash[\alpha].
\]
We need to show that if $\Gamma\vdash_{\!{\mathsf{AS}}}\Delta$, then $\Gamma\vdash\Delta$. We prove it by induction on $n=\sum\limits_{\psi \in\vert \Delta\vert }\Delta(\psi)  $. If $n=0$ the claim is trivial. Assume now that $\Delta=\Delta_{0}\uplus\lbrack\varphi]$. By Lemma~\ref{l:trees}, there are multisets $\Gamma^{\varphi}$ and $\Gamma^{r}$ such that $\Gamma^{\varphi}\uplus\Gamma^{r}=\Gamma$, $\Gamma^{\varphi} \vdash_{\!{\mathsf{AS}}}^{t}\varphi$ and $\Gamma^{r}\vdash_{\!{\mathsf{AS}} }\Delta\setminus\lbrack\varphi]$. Thus, our assumption on $\vdash$ and the induction hypothesis imply that $\Gamma^{\varphi}\vdash\varphi$ and $\Gamma^{r} \vdash\Delta\setminus[\varphi]$. Thus $\Gamma\vdash[\varphi ]\uplus\Gamma^{r}$ and $\Gamma^{r} \uplus[\varphi] \vdash\Delta$. Transitivity completes the proof.
\end{proof}

We close this subsection by providing an axiomatic system for $\vdash _{\mathcal{MV}}$, the multiset companion of infinite-valued \L{}ukasiewicz logic $\vdash_{\text{\L }}$. Although we show that it has no single\=/conclusion axiomatisation, we also axiomatise its single\=/conclusion fragment. In this way, we incidentally provide an example of two different \mdr's with the same single\=/conclusion fragment.

\begin{proposition}
There is no single\=/conclusion presentation of\/ $\vdash_{\mathcal{MV}}$.
\end{proposition}

\begin{proof}
Assume that there is such a system $\mathsf{Ax}$ and note that we have $[p\otimes q]\vdash_{\mathsf{Ax}}[p,q]$. Then, due to Lemma~\ref{l:trees}, either $\vdash_{\mathsf{Ax}}p$ or $\vdash_{\mathsf{Ax}}q$, a contradiction.
\end{proof}

\begin{definition}
The axiomatic system $\mathsf{MV}^{s}$, formulated in the language $\mathcal{L}_{0}$ of commutative residuated lattices, contains as axioms all instances of the axioms of \L{}ukasiewicz logic in $\mathcal{L}_0$, and as its sole deduction rule the rule (MP): $[\varphi,{\varphi \to\psi}] \rhd\lbrack\psi ]$. The axiomatic system $\mathsf{MV}$ is an extension of $\mathsf{MV}^{s}$ by the rule ($\otimes$-Elim): $[\varphi\otimes\psi] \rhd \lbrack\varphi,\psi]$.
\end{definition}

Observe that $\mathsf{MV}^{s}$ is a single\=/conclusion axiomatic system.

\begin{theorem}
Let $\Gamma,\Delta,[\varphi]\in \Fm_{\mathcal{L}_{0}}^{\flat}$. Then, we have the following:
\begin{enumerate}
\item $\Gamma\vdash_{\mathsf{MV}^{s}} [\varphi] \iff \Gamma\vdash_{\mathcal{MV}}[\varphi]$.

\item $\Gamma\vdash_{\mathsf{MV}}\Delta \iff \Gamma\vdash_{\mathcal{MV}}\Delta$.
\end{enumerate}
\end{theorem}

\begin{proof}
Recall, that given a multiset of $\mathcal{L}_{0}$\=/formulas $\Gamma = [\gamma _{1},\dots,\gamma_{n}]$, we write $\bigotimes\Gamma$ for $\gamma _{1}\otimes\cdots\otimes\gamma_{n}$. Note that due to the standard completeness of $\vdash_{\text{\L}}$ we obtain:
\begin{equation}
 \Gamma \vdash_{\mathcal{MV}} \Delta
 \quad\iff\quad
 {\mathcal{MV}} \models \bigotimes \Gamma \leq \bigotimes\Delta
 \quad\iff\quad\vdash_{\text{\L}} \bigotimes \Gamma\to\bigotimes\Delta.
\end{equation}
For the left-to-right direction of the former claim, it suffices to prove the latter. Assume that $\Gamma=\Gamma_{1},\dots,\Gamma_{n}\geq\Delta$ is a derivation of $\Delta$ from $\Gamma$. If we show that ${\mathcal{MV}} \models\bigotimes\Gamma_{i}\leq\bigotimes\Gamma_{i+1}$ the claim follows as ${\mathcal{MV}}\models\bigotimes\Gamma_{n}\leq\bigotimes\Delta$. We distinguish three cases:
\begin{enumerate}
\item[$\bullet$] The case when $\Gamma_{i+1}=\Gamma_{i}\uplus\lbrack\varphi]$, where $\varphi$ is an axiom, is simple, as in this case we have ${\mathcal{MV} }\models\varphi \approx 1$.

\item[$\bullet$] The case when there is a multiset $\Delta$ such that $\Gamma_{i} = \Delta\uplus[\varphi,\varphi\to\psi]$ and $\Gamma_{i+1} = \Delta\uplus[\psi]$: we know that ${\mathcal{MV}}\models\varphi\otimes(\varphi\to\psi) \leq\psi$ and so ${\mathcal{MV}}\models\varphi\otimes(\varphi\to\psi) \otimes \bigotimes\Delta\leq\psi\otimes\bigotimes\Delta$.

\item[$\bullet$] The final case, when there is a multiset $\Delta$ such that $\Gamma _{i}=\Delta\uplus\lbrack\varphi\otimes\psi]$ and $\Gamma_{i+1}=\Delta \uplus\lbrack\varphi,\psi]$, is simple.
\end{enumerate}

For the converse direction, we first prove the first claim by induction on the length of~$\Gamma$. Note that thanks to Lemma~\ref{l:trees} we can work with tree-proofs. If $\Gamma=\emptyset$, then by the assumption we know that $\varphi$ is a theorem of \L{}ukasiewicz logic, i.e., there is a proof of $\varphi$ in its usual Hilbert calculus. This proof can be easily transformed into a tree-proof of $\varphi$ in $\vdash_{\mathsf{MV} ^{s}}$. For the induction step, observe that if $\Gamma\uplus\lbrack \psi]\vdash_{\mathcal{MV}}[\varphi]$, then $\Gamma\vdash_{\mathcal{MV}} [\psi \to \varphi]$ and so by induction there is a tree-proof of $\psi \to \varphi$ from $\Gamma$ from which it is trivial to get a tree-proof of $\varphi$ from $\Gamma\uplus\lbrack\psi]$. Finally, we prove the right-to-left direction of the latter claim. First notice that $\Gamma \vdash_{\mathcal{MV}}\Delta$ entails $\Gamma\vdash_{\mathcal{MV}} \bigotimes\Delta$ and so by the first claim $\Gamma\vdash_{\mathsf{MV}^{s} }\bigotimes\Delta$ and thus also $\Gamma\vdash_{\mathsf{MV}}\bigotimes\Delta$. Repeated use of the rule ($\otimes$-Elim) completes the proof.
\end{proof}

\subsection*{Acknowledgements}

This work has received funding from the European Union's Horizon 2020 research and innovation programme under the Marie Sklodowska-Curie grant agreement No 689176 (project ``Syntax Meets Semantics: Methods, Interactions, and Connections in Substructural logics''). Petr Cintula and Tommaso Moraschini were also supported by RVO 67985807 and Czech Science Foundation project~GBP202/12/G061. Petr Cintula also acknowledges financial support from Regione Autonoma Sardegna, Visiting Professors Programme~2017. Jos\'e Gil-F\'erez also acknowledges the financial support of the Swiss National Science Foundation (SNF) grant 200021\_165850.

Conversations with Libor B\v{e}hounek and Constantine Tsinakis were extremely fruitful in shaping the contents of the paper. Preliminary versions of this work were presented in conferences in Barcelona, Bochum, Melbourne, Nashville, Storrs. We thank the audiences of these talks for providing very useful feedback. Finally, let us thank the anonymous referee for her or his helpful comments.

\end{document}